\documentclass[12pt]{amsart}
\usepackage{a4wide,graphicx}
\usepackage{epstopdf}
\allowdisplaybreaks

\let\pa\partial  
\let\na\nabla  
\let\eps\varepsilon  
\let\ds\displaystyle
\newcommand{\N}{{\mathbb N}}  
\newcommand{\R}{{\mathbb R}} 
\newcommand{\C}{{\mathbb C}} 
\newcommand{\diver}{\operatorname{div}}

\newcommand{\MM}{\mathcal{M}}
\newcommand{\TT}{\mathcal{T}}
\newcommand{\PP}{\mathcal{P}}
\newcommand{\EE}{\mathcal{E}}
\newcommand{\DD}{\mathrm{D}}
\newcommand{\dist}{\operatorname{d}}
\newcommand{\m}{\operatorname{m}}

\newtheorem{theorem}{Theorem}

\usepackage{pgfplots}           % graphics from MATLAB

\pgfplotsset{compat=newest}
\pgfplotsset{plot coordinates/math parser=false}

\newlength\entropywidth
\newlength\vacwidthneg
\newlength\vacwidthall
\newlength\cltopwidth
\newlength\optoclwidth
\newlength\cltopspwidth
\newlength\optoclspwidth
\begin{document}
\title[A finite-volume scheme for a spinorial drift-diffusion model]{A finite-volume
scheme for a spinorial \\ matrix drift-diffusion model for semiconductors}

\author[C. Chainais-Hillairet]{Claire Chainais-Hillairet}
\address{Laboratoire Paul Painlev\'e, Universit\'e Lille 1 Sciences et Technologies,
Cit\'e Scientifique, 59655 Villeneuve d'Ascq Cedex, France}
\email{Claire.Chainais@math.univ-lille1.fr}

\author[A. J\"ungel]{Ansgar J\"ungel}
\address{Institute for Analysis and Scientific Computing, Vienna University of
	Technology, Wiedner Hauptstra\ss e 8--10, 1040 Wien, Austria}
\email{juengel@tuwien.ac.at}

\author[P. Shpartko]{Polina Shpartko}
\address{Institute for Analysis and Scientific Computing, Vienna University of
	Technology, Wiedner Hauptstra\ss e 8--10, 1040 Wien, Austria}
\email{polina.shpartko@tuwien.ac.at}

\date{\today}

\thanks{The authors acknowledge partial support from  
the Austrian-French Project Amad\'ee of the Austrian Exchange Service (\"OAD).
The first and last author have been partially supported by 
the Austrian Science Fund (FWF), grants P24304, P27352, and W1245.
}

\begin{abstract}
An implicit Euler finite-volume scheme for a spinorial matrix drift-diffusion 
model for semiconductors is analyzed. The model consists of strongly coupled
parabolic equations for the electron density matrix or, alternatively, of
weakly coupled equations for the
charge and spin-vector densities, coupled to the Poisson equation
for the electric potential. The equations are solved in a bounded domain 
with mixed Dirichlet-Neumann boundary conditions. 
The charge and spin-vector fluxes are approximated 
by a Scharfetter-Gummel discretization. The main features of the numerical scheme
are the preservation of positivity and $L^\infty$ bounds
and the dissipation of the discrete free energy.
The existence of a bounded discrete solution and the monotonicity of the
discrete free energy are proved. For undoped semiconductor materials, the
numerical scheme is unconditionally stable. 
The fundamental ideas are reformulations
using spin-up and spin-down densities and certain projections of the spin-vector
density, free energy estimates, and a discrete Moser iteration.
Furthermore, numerical simulations of a simple ferromagnetic-layer
field-effect transistor in two space dimensions are presented.
\end{abstract}

% \paragraph{Keywords:}  
 
% \paragraph{AMS classification:}  

\maketitle

%%%%%%%%%%%%%%%%%%%%%%%%%%%%%%%%%%%%%%%%%%%%%%%%%%%%%%%%%%%%%%%%%%%%%%%%%%%%%%

\section{Introduction}

The exploitation of the electron spin in semiconductor devices is one of
the promising trends for future electronics. Since the electron current
can be controlled without changing the carrier concentration, this
may allow for (almost) energy-conserving and fast-switching devices and, 
more generally, for electronic devices based on new operating principles.
In the literature, several models have been proposed to describe the spin-polarized
transport in semiconductor structures \cite{FMESZ07,ZFD04}. 
Drift-diffusion approximations are widely employed \cite{PSP,ZFD02}, since
they do not require large computational resources but still describe the main
transport phenomena. In this paper, we aim to analyze a finite-volume scheme
for a spin drift-diffusion system. 
Before we explain the model equations, we sketch the state of the art in
spinorial drift-diffusion modeling. 

The existing drift-diffusion models can be classified into two main groups.
The first group is given by two-component 
drift-diffusion equations for the spin-up and spin-down densities.
One version of this model was rigorously derived from the spinor Boltzmann equation 
in the diffusion
limit with strong spin-orbit coupling (compared to the mean-free path)
\cite{ElH14}. A mathematical analysis
of the limit model was performed in \cite{Gli08}, proving the global-in-time
existence of weak solutions and their equilibration properties in two
space dimensions. In three space dimensions, the well-posedness of the stationary
system was shown in \cite{GaGl10}. A quantum correction of Bohm potential type
was derived in \cite{BaMe10}. 

The second group consists of
spin-vector drift-diffusion models in which the spin variable is a vector
quantity. Combining the charge density with the spin-vector density, we can
define the electron density matrix which solves a spinorial
matrix drift-diffusion system.
These models can be derived from the spinor Boltzmann equation by assuming a moderate
spin-orbit coupling \cite{ElH14}. Projecting the
spin-vector density in the direction of the precession vector, we recover
the two-component drift-diffusion system as a special case. 
In \cite{ElH14}, the scattering rates are supposed to be scalar quantities. 
Assuming that the scattering rates are positive definite Hermitian matrices,
a more general matrix drift-diffusion model was derived in \cite{PoNe11}.
The global existence of weak solutions to this model was shown in \cite{JNS15}.

The aim of this paper is to analyze an implicit Euler finite-volume
approximation of the spinorial matrix drift-diffusion model
of \cite{PoNe11} and to present some numerical simulations in two space dimensions.
A numerical analysis of a finite-volume scheme of the stationary two-component 
drift-diffusion equations was performed in \cite{GaGl10}. A finite-element
scheme for a spin-vector equation with given electron current density
(but coupled to the Landau-Lifshitz-Gilbert equation)
was analyzed in \cite{AHPPRS14} and simulated in \cite{ARBVHPS14}. 
However, no numerical analysis seems to be available so far for general 
spin-vector drift-diffusion models.

%%%%%%%%%%%%%%%%

\subsection{The model equations}\label{sec.model}

The spin-vector model of \cite{PoNe11}, which is analyzed in this paper,
consists of the drift-diffusion equation for the (Hermitian)
electron density matrix $N\in\C^{2\times 2}$ and the current density matrix
$J\in\C^{2\times 2}$,
\begin{align}
  & \pa_t N + \diver J + i\gamma[N,\vec m\cdot\vec\sigma]
  = \frac{1}{\tau}\left(\frac12\mbox{tr}(N)\sigma_0-N\right), \label{1.N} \\
  & J = -DP^{-1/2}(\na N+N\na V)P^{-1/2} \quad\mbox{in }\Omega,\ t>0, \label{1.J}
\end{align}
where $[A,B]=AB-BA$ is the commutator of two matrices $A$ and $B$ and 
$\Omega\subset\R^2$ is a bounded domain.
The scaled physical parameters are the strength of the effective magnetic field,
$\gamma>0$, the (normalized) direction of the precession vector
$\vec m=(m_1,m_2,m_3)\in\R^3$, the spin-flip relaxation time $\tau>0$, 
and the diffusion coefficient $D>0$. 
The precession vector plays the role of the local direction of the 
magnetization in the ferromagnet. 
%It is derived from the spin-orbit coupling in the
%band structure of the material and appears in the spin-orbit Hamiltonian. 
%(see \cite[Formula (67)]{ZFD04}).
In the analytic part of this paper, we assume that the precession vector
$\vec m$ is constant.
Furthermore, $\vec\sigma=(\sigma_1,\sigma_2,\sigma_3)$ is the triple of the 
Pauli matrices (see \cite[Formula (1)]{PoNe11}), $\sigma_0$ is the 
identity matrix in $\C^{2\times 2}$,
$\mbox{tr}(N)$ denotes the trace of the matrix $N$, and
$P=\sigma_0+p\vec m\cdot\vec\sigma$, where $p\in[0,1)$ represents
the spin polarization of the scattering rates. 
The product $\vec m\cdot\vec\sigma$ equals $m_1\sigma_1+m_2\sigma_2+m_3\sigma_3$.
The electric potential $V$ is self-consistently given by the Poisson equation
\begin{equation}\label{1.V}
  -\lambda_D^2\Delta V = \mbox{tr}(N)-C(x)\quad\mbox{in }\Omega, 
\end{equation}
where $\lambda_D>0$ is the scaled Debye length and $C(x)$ denotes
the doping profile \cite{Jue09}. 
The boundary and initial conditions are specified below.

In this paper, we investigate a scalar form of equations \eqref{1.N}-\eqref{1.J}.
For this, we develop $N$ and $J$ in the Pauli basis via
$N=\frac12 n_0\sigma_0 + \vec n\cdot\vec\sigma$ and 
$J=\frac12 j_0\sigma_0 + \vec j\cdot\vec\sigma$, where $n_0$ is the electron
charge density and $\vec n$ the spin-vector density. Setting $\vec n=(n_1,n_2,n_3)$
and $\vec j=(j_1,j_2,j_3)$ and defining $\eta=\sqrt{1-p^2}$, 
system \eqref{1.N}-\eqref{1.J} can be written
equivalently (see \cite[Remark 1]{PoNe11}) as
\begin{align}
  & \pa_t n_0 +\diver j_0	= 0,\label{1.n0} \\
  & \pa_t n_\ell + \diver j_\ell
	- 2\gamma(\vec n\times\vec m)_\ell = -\frac{n_\ell}{\tau},\quad 
	\ell=1,2,3, \label{1.vecn} \\
  & j_0 = \frac{D}{\eta^2}(J_0 - 2p\vec J\cdot\vec m),\quad 
	j_\ell = \frac{D}{\eta^2}\left(\eta J_\ell 
	+ (1-\eta)(\vec J\cdot\vec m)m_\ell - \frac{p}{2}J_0m_\ell\right),\quad 
	\ell=1,2,3, \label{1.vecj}\\
  & J_0 = -\na n_0-n_0\na V, \quad 
  \vec J = (J_1,J_2,J_3) = -\na\vec n - \vec n\na V \quad \mbox{in } \Omega,\ t>0.
	\label{1.J0}
\end{align}
Moreover, the Poisson equation \eqref{1.V} rewrites 
\begin{equation}\label{1.V2}
  -\lambda_D^2\Delta V = n_0-C(x)\quad\mbox{in }\Omega.
\end{equation}
System \eqref{1.n0}-\eqref{1.V2} is strongly coupled due to the cross-diffusion
terms in   \eqref{1.vecj} and nonlinear due to the Poisson coupling.
Note that any solution $(n_0,\vec n)$ to \eqref{1.n0}-\eqref{1.J0}
defines a solution $N$ to \eqref{1.N}-\eqref{1.J} and vice versa.

The boundary $\pa\Omega=\Gamma^D\cup\Gamma^N$ is assumed to consist of
the union of contacts $\Gamma^D$ and the isolating boundary part $\Gamma^N$.
Then the  boundary and initial data are given by
\begin{align}
  & n_0=n^D, \quad \vec n=0, \quad V=V^D\quad \mbox{on }\Gamma^D,\ t>0, 
	\label{1.bcD} \\
	& \na n_0\cdot\nu =  \na n_\ell\cdot\nu=\na V\cdot\nu = 0
	\quad\mbox{on }\Gamma^N,\ t>0, \ \ell=1,2,3, \label{1.bcN} \\
  & n_0(0)=n_0^0, \quad \vec n(0)=\vec{n}^0\quad\mbox{in }\Omega, \label{1.ic}
\end{align}
where $\nu$ is the exterior unit normal vector to $\pa\Omega$.

%%%%%%%%%%%%%%%%%

\subsection{Mathematical background}\label{sec.back}

We aim to design a numerical scheme which preserves some qualitative
properties of the continuous model, in particular preservation of the 
positivity of the charge density, boundedness of the density matrix,
and dissipation of the free energy. The main difficulty of the analysis 
is the strong coupling of the equations (the diffusion matrix is not diagonal), 
since maximum principle or regularity arguments generally do not apply.
The key idea is to introduce two transformations of variables which
make the diffusion matrix diagonal and thus reduce the level of coupling.

The first transformation is defined by the spin-up and spin-down densities
$n_\pm = \frac12 n_0\pm\vec n\cdot\vec m$. Then system \eqref{1.n0}-\eqref{1.J0}
becomes
\begin{align}
 & \pa_t n_+ + \diver\big(D(1+ p)(-\na n_+ - n_+\na V)\big) 
	= -\frac{1}{2\tau}(n_+ - n_-), \label{1.np}\\
 & \pa_t n_- + \diver\big(D(1- p)(-\na n_- - n_-\na V)\big) 
	= -\frac{1}{2\tau}(n_- - n_+) \label{1.nm}	
\end{align}
and the boundary conditions \eqref{1.bcD}, \eqref{1.bcN} imply 
\begin{equation}
  n_\pm=\frac{n^D}{2} \mbox{ on }\Gamma^D\quad\mbox{and}\quad 
	\nabla n_\pm \cdot \nu=0 \mbox{ on }\Gamma^N,\ t>0.\label{1.bc.npm}
\end{equation}

We observe that \eqref{1.n0}-\eqref{1.J0} implies \eqref{1.np}-\eqref{1.nm} 
but not vice versa. 
Physically this is clear since the spin-up and spin-down densities contain less 
information than the full density matrix $N$. 
By the Stampacchia truncation method, the positivity and boundedness of 
$n_\pm$ was shown in \cite{JNS15}, thus giving the positivity and boundedness
of the charge density $n_0=n_++n_-$. 
Using the notation $\sum_\pm a_\pm=a_+ + a_-$, the (relative) free energy 
of the above system is given by the sum of the entropy and the electric energy,
\begin{equation}\label{1.E}
  E(t) = \sum_\pm\int_\Omega\left(n_\pm(\log n_\pm-1) - n_\pm\log\frac{n^D}{2} 
	+ \frac{n^D}{2}\right)dx 
	+ \frac{\lambda_D^2}{2}\int_\Omega|\na(V-V^D)|^2 dx.
\end{equation}
Some formal computations show that it is nonnegative and nonincreasing for $t>0$.

The second transformation is given by the decomposition of $\vec n$ in the
parallel and perpendicular components with respect to $\vec m$: 
$\vec n_\parallel=(\vec n\cdot\vec m)\vec m$ and 
$\vec n_\perp=\vec n-\vec n_\parallel$. The equation for $\vec n_\perp$ reads as
\begin{equation}\label{1.nperp}
  \pa_t \vec n_\perp +\diver\left(\frac{D}{\eta}(-
  \na\vec n_\perp - \vec n_\perp\na V)\right) - 2\gamma(\vec n_\perp\times\vec m) 
	= - \frac{\vec n_\perp}{\tau}. 
\end{equation}
In \cite{JNS15}, it was proved by a Moser-type iteration technique that
$\vec n_\perp$ is bounded. Since $\vec n_\parallel=\frac12(n_+-n_-)\vec m$ is
bounded as well (see above), this implies an $L^\infty$ bound for $\vec n$
and consequently for the density matrix $N=\frac12 n_0\sigma_0+\vec n\cdot\vec\sigma$.

The task is to ``translate'' these ideas to a finite-volume setting.
We approximate the diffusive and convective part of the fluxes
simultaneously by using a Scharfetter-Gummel discretization.
These fluxes were introduced by Il'in \cite{Ili69} and Scharfetter and Gummel
\cite{ScGu69} for the classical drift-diffusion model (without spin coupling). 
The discretizations are second-order accurate in space and preserve the
steady states. The dissipativity with an implicit Euler discretization
was shown in \cite{GaGa96}. The discrete steady states were proved to be bounded
\cite{GaGl10}. Discrete entropy (free energy) estimates and/or the exponential decay
of the free energy along trajectories towards the global equilibrium were 
investigated in \cite{Cha11,Gli11} but still without any spin coupling.

Our main results, detailed in Section \ref{sec.nummain}, are the existence
of a bounded discrete solution to a fully discrete finite-volume scheme for  
\eqref{1.n0}-\eqref{1.ic} and the monotonicity
of the discrete free energy for the spin-up and spin-down densities.
The mathematical challenge is the proof of lower and upper bounds for the 
discrete densities.
The ``translation'' of the Stampacchia truncation argument from the continuous
to the discret case in, e.g., \eqref{1.nperp} faces some
difficulties due to the drift term. The main difficulty lies in the fact
that the monotonicity of the drift term (with respect to the density
variable) cannot be exploited. Therefore, a Moser-type
iteration method was employed in \cite{JNS15}. The idea is to derive a uniform
estimate for $\vec n_\perp$ in the $L^q$ norm of the form
$$
  \frac{d}{dt}\|\vec n_\perp\|_q^q \le cq\|\vec n_\perp\|_q^q, \quad t>0,
$$
where $\|\cdot\|_q$ denotes the $L^q(\Omega)$ norm and $c>0$ does not depend on
$q\in(1,\infty)$. By Gronwall's lemma, this implies that
$$
  \|\vec n(t)\|_q \le e^{ct}\|\vec n_\perp(0)\|_q,
$$
and the limit $q\to\infty$ shows the claim. The discrete
equivalent of the estimate is
$$
  \|\vec n_\perp^k\|_q^q - \|\vec n_\perp^{k-1}\|_q^q
	\le cq\triangle t\|\vec n_\perp^k\|_q^q,
$$
where $\vec n_\perp^k$ is an approximation of $\vec n_\perp$ at time $t^k$
and $\triangle t$ is the (uniform) time step size.
In order to solve this recursion, we require $1-cq\triangle t>0$,
thus imposing a condition on the time step size for fixed $q$. This motivates
additional conditions on the model parameters, which are described and
discussed in Section \ref{sec.main}.

The paper is organized as follows. In Section \ref{sec.nummain}, we detail the
numerical scheme and present the main results, in particular the existence
of discrete solutions (Theorem \ref{thm.ex}) and the dissipativity of the
discrete free energy (Theorem \ref{thm.diss}). The proofs are given in
Sections \ref{sec.proof1} and \ref{sec.proof2}. Some numerical tests are presented
in Section \ref{sec.numer}.

%%%%%%%%%%%%%%%%%%%%%%%%%%%%%%%%%%%%%%%%%%%%%%%%%%%%%%%%%%%%%%%%%%%%%%%%%%%%%%%

\section{Numerical method and main results}\label{sec.nummain}

In this section, we specify the numerical discretization of the spin drift-diffusion
system \eqref{1.n0}-\eqref{1.ic} and state the main results of the paper.

\subsection{Notations}\label{sec.not}

Before we state the numerical scheme, we need to define the mesh of the
domain $\Omega$ and to introduce some notation. 
We consider the two-dimensional case only but the scheme can
be generalized in a straightforward way to higher dimensions. 

Let $\Omega\subset\R^2$ be an open bounded polygonal set. The mesh
$\MM=(\TT,\EE,\PP)$ is given by a family $\TT$ of open polygonal control volumes
or cells, a family $\EE$ of edges, and a family $\PP=(x_K)_{K\in\TT}$ of points.
We assume that the mesh is admissible in the sense of Definition 9.1 in \cite{EGH00}.
This definition implies that the straight line between two neighboring centers
of cell $(x_K,x_L)$ is orthogonal to the edge $\sigma=K|L$ between two control
volumes $K$ and $L$ and therefore collinear to the unit normal vector 
$\nu_{K,\sigma}$ to $\sigma$
outward to $K$. For instance, triangular meshes satisfy the admissibility condition
if all angles of the triangles are smaller than $\pi/2$ \cite[Example 9.1]{EGH00}.
Voronoi meshes are also admissible meshes \cite[Example 9.2]{EGH00}.

Each edge $\sigma\in\EE$ is either an internal edge, $\sigma=K|L$, or an exterior
edge, $\sigma\subset\pa\Omega$, and we set $\EE=\EE_{\rm int}\cup\EE_{\rm ext}$. 
We assume that each exterior edge is an element of either the Dirichlet or
Neumann boundary such that we can set 
$\EE_{\rm ext}=\EE_{\rm ext}^D\cup\EE_{\rm ext}^N$.
For a given control volume $K\in\TT$, we define the set $\EE_K$ of the edges of $K$,
which can be written as the union of $\EE_{K,{\rm int}}$, $\EE_{K,{\rm ext}}^D$,
and $\EE_{K,{\rm ext}}^N$. For every $\sigma\in\EE$, there exists at least one cell
$K\in\TT$ satisfying $\sigma\in\EE_K$, and we denote this cell by $K_\sigma$. 
When $\sigma$ is an interior cell with $\sigma=K|L$, we have $K_\sigma=K$ or
$K_\sigma=L$.

For $K\in\TT$ and $\sigma\in \EE_K$, we denote by $\dist_{K,\sigma}$ the 
distance $\dist_{K,\sigma}=d(x_K,\sigma)$. Then, for $\sigma\in\EE_{\rm int}$, 
$\sigma=K\vert L$, we define $\dist_\sigma=\dist_{K,\sigma}+\dist_{L,\sigma}
=d(x_K,x_L)$ and for $\sigma\in \EE_{\rm ext}$ 
with $\sigma\in \EE_K$, $\dist_\sigma=\dist_{K,\sigma}$. 
Furthermore, the measure of $\sigma\in\EE$ or a set 
$\omega\subset\Omega$ is denoted by $\m(\sigma)$ or $\m(\omega)$, respectively.
In the numerical scheme, we need the so-called transmissibility coefficient
$\tau_\sigma=\m(\sigma)/\dist_\sigma$ for $\sigma\in\EE$. We assume that the
mesh satisfies the regularity constraint
\begin{equation}\label{sch.mesh}
  \exists\xi>0:\ \forall K\in\TT:\ \forall \sigma\in\EE_K:\ 
	\dist(x_K,\sigma)\ge \xi\mbox{diam}(K).
\end{equation}

The finite-volume scheme for a conservation law with unknown $u$ provides a
vector $u_\TT=(u_K)_{K\in\TT}$ of approximate values and the associated piecewise
constant function, still denoted by $u_\TT$, $u_\TT=\sum_{K\in\TT}u_K\mathbf{1}_K$,
which approximates the unknown $u$. Here, $\mathbf{1}_K$ denotes the 
characteristic function of the cell $K$. The approximate values of the Dirichlet 
boundary provide a vector $u_{\EE^D}=(u_\sigma)_{\sigma\in\EE_{\rm ext}^D}$.
The vector containing the approximate values in the control volumes and at the
Dirichlet boundary edges is denoted by $u_\MM=(u_\TT,u_{\EE^D})$. 

The numerical scheme can be formulated in a compact form by introducing the
following notation. For any vector $u_\MM=(u_\TT,u_{\EE^D})$, we define,
for all $K\in\TT$ and $\sigma\in\EE_K$,
$$
  u_{K,\sigma} = \left\{
  \begin{aligned}
    u_L      &\quad\mbox{if } \sigma = K | L, \\
    u_\sigma &\quad\mbox{if } \sigma = \EE_{K,{\rm ext}}^D, \\
    u_K      &\quad\mbox{if } \sigma = \EE_{K,{\rm ext}}^N,
  \end{aligned}\right.
$$
and we set $ \DD u_{K,\sigma}=u_{K,\sigma}-u_K$. We remark that the definition 
of $u_{K,\sigma}$ ensures that $\DD u_{K,\sigma}=0$ on the Neumann boundary edges. 
Then the discrete $H^1$ seminorm for $u_\MM$ can be defined by
$$
  |u_\MM|_{1,\MM} = \left(\sum_{\sigma\in\EE}\tau_\sigma
	|\DD u_{K,\sigma}|^2\right)^{1/2},
$$
where the summation is over all edges $\sigma\in\EE$ with $K=K_\sigma$.
The $L^p$ norm of $u_\TT$ reads as
$$
  \|u_\TT\|_p = \left(\sum_{K\in\TT}\m(K)|u_K|^p\right)^{1/p}\mbox{ for } 
	1\le p<\infty \quad\mbox{and}\quad  \|u_\TT\|_\infty=\ds\max_{K\in\TT} |u_K|.
$$

When formulating a finite-volume scheme, we have to define some numerical fluxes 
$J_{K,\sigma}$ which are consistent approximations of the exact fluxes 
through the edges $\int_\sigma J\cdot \nu_{K,\sigma}ds$. 
We impose the conservation of the numerical fluxes $J_{K,\sigma}+J_{L,\sigma}$ 
for $\sigma =K|L$, requiring that they 
vanish on the Neumann boundary edges, $J_{K,\sigma}=0$ for 
$\sigma \in \EE_{K,{\rm ext}}^N$. 
Then the discrete integration-by-parts formula becomes
\begin{equation}\label{ibp}
  \sum_{K\in\TT}\sum_{\sigma\in\EE_K} J_{K,\sigma}u_K
	= -\sum_{\sigma\in\EE} J_{K,\sigma}\DD u_{K,\sigma} 
	+ \sum_{\sigma\in\EE_{\rm ext}^D} J_{K,\sigma} u_{K,\sigma}.
\end{equation}

%%%%%%%%%%%%%%

\subsection{Numerical scheme}\label{sec.scheme}

At each time step $k\ge 0$, we define the approximate solution
$u^k_\TT=(u^k_K)_{K\in\TT}$ for $u\in\{n_0,\vec{n},V\}$ and the approximate values 
at the Dirichlet boundary, $u^k_{\EE^D}=(u^k_\sigma)_{\sigma\in\EE_{\rm ext}^D}$
(which in fact does not depend on $k$ since the boundary data is time-independent).
We first define the initial and boundary conditions corresponding to 
\eqref{1.ic} and \eqref{1.bcD}. We set
\begin{align}
  (n_{0,K}^0,\vec{n}_K^0) &= \frac{1}{\m(K)}\int_K(n_0^0,\vec{n}^0)dx
	\quad\mbox{for all }K\in\TT, \label{sch.ic}\\
	(n_{0,\sigma}^D,\vec{n}_\sigma^D,V_\sigma^D) 
	&= \frac{1}{\m(\sigma)}\int_\sigma(n^D,\vec{0}, V^D)ds \quad\mbox{for all }
	\sigma\in\EE_{\rm ext}^D.\nonumber
\end{align}
Note that 
$\vec{n}_\sigma^D=0$ for $\sigma\in\EE_{\rm ext}^D$. We may define similarly 
the quantities $\vec{m}_K, C_K, D_K, p_K$ for a given $K\in\TT$.

We consider a temporal implicit Euler and spatial finite-volume discretization.
The scheme for \eqref{1.n0}, \eqref{1.vecn}, \eqref{1.V2} writes, 
for all $K\in\TT$ and $k\geq 1$, as
\begin{align}
  & \m(K)\frac{n_{0,K}^k-n_{0,K}^{k-1}}{\triangle t} 
	+\sum_{\sigma\in\EE_K}j_{0,K,\sigma}^k = 0, \label{sch.n0} \\
	& \m(K)\frac{\vec{n}_K^k-\vec{n}_K^{k-1}}{\triangle t} 
	+ \sum_{\sigma\in\EE_K}\vec{j}_{K,\sigma}^k 
	- 2\gamma\m(K)(\vec{n}_K^k\times\vec{m}_K) = -\frac{\m(K)}{\tau}\vec{n}_K^k, 
	\label{sch.nl} \\
	&-\lambda_D^2\sum_{\sigma\in\EE_K}\tau_\sigma \DD V_{K,\sigma}^k
	= \m(K)(n_{0,K}^k-C_K), \label{sch.v}
\end{align}
where the discrete counterpart to \eqref{1.vecj} is, for all $K\in\TT$, 
$\sigma\in\EE_K$, $k\geq 0$,
\begin{align}
  j_{0,K,\sigma}^k &= \frac{D_\sigma}{\eta_\sigma^2}\big(J_{0,K,\sigma}^k
	-2p_\sigma \vec J_{K,\sigma}^k\cdot\vec{m}_\sigma\big), \label{sch.j0} \\
	\vec{j}_{K,\sigma}^k &= \frac{D_\sigma}{\eta_\sigma^2}\Big(
	\eta_\sigma \vec{J}_{K,\sigma}^k + (1-\eta_\sigma)(\vec{J}_{K,\sigma}^k\cdot
	\vec{m}_\sigma)\vec{m}_\sigma - \frac{p_\sigma}{2}J_{0,K,\sigma}^k
	\vec{m}_\sigma\Big). \label{sch.jl}
\end{align}
The numerical fluxes $J_{0,K,\sigma}^k$ and $J_{\ell,K,\sigma}^k$ are
approximations of the integrals $\int_\sigma J_0\cdot\nu_{K,\sigma}ds$
and $\int_\sigma J_\ell\cdot\nu_{K,\sigma}ds$ at time $k\triangle t$, 
and we set $\vec{J}_{K,\sigma}^k=(J_{\ell,K,\sigma}^k)_{\ell=1,2,3}$.
We recall that $J_0$ and $\vec{J}$ are defined by \eqref{1.J0}.
We use a Scharfetter-Gummel approximation for the definition of the numerical fluxes.
For given $K\in\TT$ and $\sigma\in\EE_K$, we set
\begin{equation}
  J_{\ell,K,\sigma}^k = \tau_\sigma\big(B(\DD V_{K,\sigma}^k)n_{\ell,K}^k
	- B(-\DD V_{K,\sigma}^k)n_{\ell,K,\sigma}^k\big), \quad \ell=0,1,2,3,
	\label{sch.J}
\end{equation}
where $B$ is the Bernoulli function defined by 
$$
  B(x)=\ds\frac{x}{\exp(x)-1} \mbox{ for }x\neq 0
  \quad\mbox{and}\quad B(0)=1.
$$ 

It remains to define the quantities $D_\sigma$, $\vec{m}_\sigma$, $p_\sigma$ and
$\eta_\sigma$ appearing in \eqref{sch.j0} and \eqref{sch.jl}.  
We use a weighted harmonic average 
on the interior edges and a classical mean value on the boundary edges,
$$
  D_\sigma 	= \frac{\dist_{\sigma} D_K D_L}{\dist_{K,\sigma} D_L+\dist_{L,\sigma} D_K}
  \mbox{ for }\sigma \in \EE_{{\rm int}},\ \sigma=K|L, \quad
  D_\sigma=\frac{1}{\m(\sigma)}\int_{\sigma} D(s) ds\mbox{ for }
	\sigma \in \EE_{{\rm ext}}^D,
$$
and similar definitions for $\vec{m}_\sigma$ and $p_\sigma$. Furthermore,
we set $\eta_{\sigma}=\sqrt{1-p_\sigma^2}$.

Finally, the boundary conditions are
\begin{align}
  n_{0,\sigma}^k = n_{0,\sigma}^D, \quad \vec{n}_\sigma^k = 0,
	\quad V_\sigma^k = V_\sigma^D &\quad\mbox{for }\sigma\in\EE_{\rm ext}^D, 
	\label{sch.bc0} \\
	\DD n_{\ell,K,\sigma}^k = \DD V_{K,\sigma}^k = 0 &\quad\mbox{for }
	\sigma\in\EE_{K,{\rm ext}}^N,\ \ell=0,1,2,3,\ k\ge 0. \label{sch.bcl}
\end{align}
We remark that they imply $J_{l,K,\sigma}^k=0$ for $\sigma\in\EE_{K,{\rm ext}}^N$, 
$\ell=0,1,2,3$, and $ k\ge 0$.

For later use, we note that, using the elementary property $B(x)-B(-x)=-x$ for 
$x\in\R$, the numerical fluxes can be reformulated in two different manners: 
\begin{align}
  J_{\ell,K,\sigma}^k &= \tau_\sigma \big(-\DD V_{K,\sigma}^k n_{\ell,K}^k
	- B(-\DD V_{K,\sigma}^k)\DD n_{\ell,K,\sigma}^k\big)  \label{1.refJ.1}\\
	&= \tau_\sigma \big(-\DD V_{K,\sigma}^k n_{\ell,K,\sigma}^k
	- B(\DD V_{K,\sigma}^k)\DD n_{\ell,K,\sigma}^k\big), \quad 
	\ell=0,1,2,3,\label{1.refJ.2}
\end{align}
and adding these expressions leads to a third formulation:
%\begin{equation}\label{sch.coth}
%  J_{\ell,K,\sigma}^k = \tau_\sigma \left(- \frac12(n_{\ell,K}^k+n_{\ell,K,\sigma}^k)\DD V_{K,\sigma}^k-\frac12 \DD V_{K,\sigma}^k
%	B^s(\DD V_{K,\sigma}^k)\DD n_{\ell,K,\sigma}^k\right),
%\end{equation}
\begin{equation}\label{sch.coth}
  J_{\ell,K,\sigma}^k = \tau_\sigma \left(- \frac12(n_{\ell,K}^k+n_{\ell,K,\sigma}^k)\DD V_{K,\sigma}^k-	B^s(\DD V_{K,\sigma}^k)\DD n_{\ell,K,\sigma}^k\right),
\end{equation}
where
\begin{equation}\label{Bs}
  B^s(x)= \frac{x}{2}\coth\left(\frac{x}{2}\right)=\frac{B(x)+B(-x)}{2}.
\end{equation}

%%%%%%%%%%%%%%%%%%%%%%%%

\subsection{Main results}\label{sec.main}

We impose the following assumptions on the domain and the data:
\begin{align}
  & \Omega\subset\R^2\mbox{ bounded domain},\ \pa\Omega=\Gamma^D\cup\Gamma^N,
	\ \Gamma^D\cap\Gamma^N=\emptyset, \ \m(\Gamma^D)>0,\ \Gamma^N
	\mbox{ open}, \label{hypo1} \\
 	&  D,\,p,\,\vec m\mbox{ are constant and } |\vec m|=1,\label{hypo2}\\
  &  n_0^0,\ \vec n^0,\ n^D\in L^\infty(\Omega),\  
  \frac12 n_0^0\pm {\vec n}^0\cdot m \geq 0,\ 
	n^D\geq 0,\ n^D,\ V^D\in H^1(\Omega), \label{hypo3} \\
	& \MM=(\TT,\EE,\PP) \mbox{ is an admissible mesh satisfying \eqref{sch.mesh}}.	
	\label{hypo4}
\end{align}

We first remark that if $(n^k_{0,\TT},\vec n^k_\TT,V^k_\TT)$ is a solution to 
scheme \eqref{sch.n0}-\eqref{sch.bcl}
for a given $k\geq 1$ ($(n^0_{0,\TT},\vec n^0_\TT)$ are defined as the 
discretization of the initial conditions), we can define
$n^k_{\pm,\TT}=\frac12 n^k_{0,\TT}\pm \vec n^k_\TT\cdot\vec m$, 
${\vec n}_{\perp,\TT}^k= \vec n^k_\TT-(\vec n^k_\TT\cdot\vec m)\vec m$.
Moreover, as $n^D$ and $V^D$ are defined on the whole domain $\Omega$, we can 
define $n_\TT^D$ and $V_\TT^D$ by taking the mean value of $n^D$ and $V^D$ on each 
control volume $K\in\TT$.

Then the following existence result holds.

\begin{theorem}[Existence of a solution to the numerical scheme and $L^\infty$ 
bounds]\label{thm.ex}
{\sloppy Let assumptions \eqref{hypo1}-\eqref{hypo4} hold.} 
We impose the following constraints:
\begin{equation}\label{hypo5}
  0 < \triangle t \le \frac{1}{\alpha} := \frac{\lambda_D^2}{D(1+p)\|C\|_\infty}, \quad
	0 < \tau \le \frac{\eta\lambda_D^2}{D\|C\|_\infty}.
\end{equation}
Then for $k\geq 1$,  there exists a solution $(n^k_{0,\TT},\vec n^k_\TT,V^k_\TT)$ 
to scheme \eqref{sch.n0}-\eqref{sch.bcl} satisyfing 
$$
  0\le n^k_{0,\TT}\le 2M^0, \quad 0\le n_{\pm,\TT}^k\le M^0, \quad
	|\vec n_\TT^k|\le 2M^k\quad\mbox{in }\Omega,
$$
where $M^k = M^0(1-\alpha\triangle t)^{-k}$ and
$$
  M^0 = \ds\max\left(\frac12\sup_{\pa \Omega} n^D,\ 
	\sup_{\Omega}\left(\frac12 n_0^0+|{\vec n}^0\cdot {\vec m}|\right),\
	\sup |{\vec n}_\perp^0|,\ \sup_{\Omega} C  \right).
$$
\end{theorem}

In the continuous case, similar $L^\infty$ bounds for the spin-up and 
spin-down densities, and therefore for the electron charge density, 
were shown in \cite{JNS15}. These bounds do not depend on time.
The mixing of the spin-vector components prevents the use
of the monotonicity argument for $\vec n_\perp$, solving \eqref{1.nperp}.
Therefore, both in the continuous and discrete situations, the $L^\infty$
bound for the spin-vector density depends on time.

The constraint on $\triangle t$ is needed in the definition of $M^k$.
%, leading to an upper bound for $n_{\pm\TT}^k$
 Furthermore, the condition on $\tau$ is
necessary to prove the $L^\infty$ bound for $\vec n_{\perp,\TT}^k$.
We believe that the latter restriction is technical. We stress the fact
that our scheme is unconditionally stable if the semiconductor is undoped,
i.e.\ $C=0$. In this situation, $\triangle t$ and $\tau$ can be chosen arbitrarily.

Next, we prove that the scheme dissipates the discrete free energy, defined by
\begin{align}
  E^k &= \sum_\pm\sum_{K\in\TT}\m(K)\left(n_{\pm,K}^k (\log n_{\pm,K}^k-1)
	- n_{\pm,K}^k\log\frac{n_{K}^D}{2} +\frac{n_{K}^D}{2}\right) 
	\nonumber \\
	&\phantom{xx}{}
	+ \frac{\lambda_D^2}{2}\sum_{\sigma\in\EE}\tau_\sigma(\DD (V^k-V^D)_{K,\sigma})^2.
	\label{sch.E}
\end{align}

\begin{theorem}[Dissipation of the discrete free energy]\label{thm.diss}
Let assumptions \eqref{hypo1}-\eqref{hypo4} hold
and let $(n_{0,\TT}^k,\vec n_\TT^k,V_\TT^k)_{k\geq 0}$ be a solution to scheme
\eqref{sch.n0}-\eqref{sch.bcl} satisyfing $0\le n_{\pm,\TT}^k\le M^0$. We further 
assume that $n^D\geq n_\ast>0$ and that $\log(n^D/2)+V^D$ is constant in 
$\overline\Omega$. Then the mapping $k\mapsto E^k$ is nonincreasing, i.e.,
the scheme dissipates the free energy \eqref{sch.E}:
\begin{equation}\label{ineq.thm.diss}
  E^k + \frac{\triangle t}{2}\sum_\pm D(1\pm p)\sum_{\sigma\in\EE}\tau_\sigma
	\min\{n^k_{\pm,K},n^k_{\pm,K,\sigma}\}
	\big(\DD(\log n^k_{\pm}+V^k)_{K,\sigma}\big)^2 \le E^{k-1},\quad k\geq 1.
\end{equation}
\end{theorem}

The above dissipation inequality for the free energy is the discrete counterpart of the
continuous estimate for the free energy \eqref{1.E} \cite[Formula (28)]{JNS15}:
$$
  E(t) + \frac12\int_0^t\int_\Omega\sum_\pm D(1\pm p)
	n_\pm|\na(\log n_\pm+V)|^2 dxds \le E(0).
$$
One may ask if the discrete solution converges to the continuous one when
the approximation parameters tend to zero. However, 
it seems to be difficult to extract a discrete gradient estimate for $n_{\pm,\TT}^k$
from the discrete free energy estimate in Theorem \ref{thm.diss} since we
do not have a suitable discrete version of the chain rule $n_\pm|\na\log n_\pm|^2
= 4|\na\sqrt{n_\pm}|^2$. 

%%%%%%%%%%%%%%%%%%%%%%%%%%%%%%%%%%%%%%%%%%%%%%%%%%%%%%%%%%%%%%%%%%%%%%%%%%%%%%%

\section{Proof of Theorem \ref{thm.ex}}\label{sec.proof1}

The proof of Theorem \ref{thm.ex} will be presented in two subsections. 
We first establish  the existence of a solution 
$(n_{0,\TT}^k, {\vec n}_\TT^k, V_\TT^k)$ at each time step $k\geq 1$ 
by an induction argument. The proof is based on the 
fixed-point theorem of Brouwer. In this subsection, we also show $L^\infty$ bounds 
on  $n_{0,\TT}^k$ and $n_{\pm,\TT}^k$ which depend on $k$. 
Then, in the second subsection, we prove that these bounds are in fact uniform with 
respect to $k$.

\subsection{Existence of a solution to the scheme}

We first note that the initial condition $(n_{0,\TT}^0, {\vec n}_\TT^0)$ is 
well-defined by \eqref{sch.ic}. Moreover, the definition of $M^0$ ensures that
$|n_{\perp,\TT}^0|\leq M^0$, $0\leq n_{\pm,\TT}^0\leq M^0$
and therefore $0\leq n_{0,\TT}^0\leq 2M^0$ and $|{\vec n}_{\pm,\TT}^0|\leq 2M^0$.

The proof is done by induction. Let $k\geq 1$. Assuming that $(n_{0,\TT}^{k-1}, 
{\vec n}_\TT^{k-1},V_\TT^{k-1})$ is given and verifies 
$|n_{\perp,\TT}^{k-1}|\leq M^{k-1}$, $0\leq n_{\pm,\TT}^{k-1}\leq M^{k-1}$, 
we will prove the existence of $(n_{0,\TT}^{k}, {\vec n}_\TT^{k}, V_\TT^{k})$, 
solution to \eqref{sch.n0}-\eqref{sch.bcl}, satisfying these bounds with $k$ 
instead of $k-1$. 
Scheme \eqref{sch.n0}-\eqref{sch.bcl} is a nonlinear system of equations. We prove the 
existence of a solution by using a fixed-point theorem. Let us denote by $\theta$ 
the cardinality of the mesh $\TT$ (the number of control volumes) and let $\mu>0$. 
We define an application $F_\mu^k: \R^{4\theta}\to \R^{4\theta}$ such that 
$F_\mu^k(\rho_\TT)=n_\TT$, where $\rho_\TT=(\rho_{0,\TT},\vec\rho_\TT)$ and 
$n_\TT=(n_{0,\TT},\vec n_\TT)$. It is based on a linearization of the scheme 
and defined in two steps:
\begin{itemize}
\item First, we define $V_\TT\in\R^{\theta}$ as the solution to  the linear system
\begin{align}
  & -\lambda_D^2\sum_{\sigma\in\EE_K}\tau_\sigma\DD V_{K,\sigma} 
	= \m(K)(\rho_{0,K}-C_K)\quad\mbox{for } K\in\TT, \label{ex.v} \\
	& V_\sigma=V_\sigma^D\mbox{ for }\sigma\in\EE_{\rm ext}^D, \quad
	\DD V_{K,\sigma}=0\mbox{ for }\sigma\in\EE_{K,{\rm ext}}^N. \nonumber
\end{align}
\item Second, we construct $n_\TT=(n_{0,\TT},\vec{n}_\TT)\in \R^{4\theta}$ 
as the solution to
\begin{align}
  \frac{\m(K)}{\triangle t}(n_{0,K} - n_{0,K}^{k-1})
	&+ \mu\frac{\m(K)}{\triangle t}(n_{0,K}-\rho_{0,K})
	+ \sum_{\sigma\in\EE_K}j_{0,K,\sigma} = 0 \quad\mbox{for } K\in\TT,\label{ex.n0} \\
	\frac{\m(K)}{\triangle t}(n_{\ell,K} - n_{\ell,K}^{k-1})
	&+ \mu\frac{\m(K)}{\triangle t}(n_{\ell,K}-\rho_{\ell,K})
	+ \sum_{\sigma\in\EE_K}j_{\ell,K,\sigma} \label{ex.nl} \\
	&{}- 2\gamma\m(K)(\vec{n}_K\times\vec{m})_\ell 
  = -\frac{\m(K)}{\tau}n_{\ell,K} \quad\mbox{for } K\in\TT,\ \ell=1,2,3, \nonumber
\end{align}
where $j_{0,K,\sigma}$ and $j_{\ell,K,\sigma}$ are defined in 
\eqref{sch.j0} and \eqref{sch.jl}, with $J_{\ell,K,\sigma}$ defined in
\eqref{sch.J}, but without the superindex $k$. The boundary conditions read as
\begin{align}
  n_{0,\sigma} = n_{\sigma}^D, \ n_{\ell,\sigma}=0
	&\quad\mbox{for }\sigma\in\EE_{\rm ext}^D,\ \ell=1,2,3, \label{ex.bc1} \\ 
	\DD n_{\ell,K,\sigma} = 0
	&\quad\mbox{for }\sigma\in\EE_{K,{\rm ext}}^N,\ K\in\TT,\
	\ell=0,1,2,3. \label{ex.bc2}
\end{align}
\end{itemize}

The parameter $\mu>0$ allows us to prove unconditional stability
for the linearized problem; see e.g.\ \cite{BCV14}. 
The corresponding term vanishes for fixed points $\rho_\TT=n_\TT$, 
so that a fixed point for $F_\mu^k$ is a solution to scheme 
\eqref{sch.n0}-\eqref{sch.bcl}. We choose 
\begin{equation}\label{ex.mu}
  \mu \ge \frac{D\|C\|_\infty}{\lambda_D^2}\max
	\left\{\frac{1}{\eta^2},\frac{1+p}{2}\right\}\triangle t.
\end{equation}

The existence and uniqueness of $V_\TT$, solution to \eqref{ex.v}, 
are obvious since the
corresponding matrix is positive definite. As this matrix does not depend on 
$\rho_\TT$ and the right-hand side is continuous with respect to $\rho_\TT$, 
the first mapping $\rho_\TT\mapsto V_\TT$ is continuous from 
$\R^{4\theta}$ to $\R^{\theta}$. This property 
is not so obvious for the second mapping, 
based on the linear system of equations \eqref{ex.n0}-\eqref{ex.bc2}. We will 
prove this property below (Step 1), in order to guarantee that the mapping 
$F_{\mu}^k$ is well-defined and continuous.

Then, in order to apply Brouwer's fixed-point theorem, we will prove that 
$F_{\mu}^k$ preserves the set 
\begin{equation}\label{def.Sk}
  {\mathcal S}^k = \left\{ n_\TT=(n_{0,\TT},{\vec n}_\TT)\in \R^{4\theta} 
	:\ 0\leq n_{\pm,\TT}\leq M^k,\ |{\vec n}_{\perp,\TT}|\leq M^k\right\}.
\end{equation}
It is a bounded set because each element $n_\TT\in {\mathcal S}^k$ verifies 
$0\leq n_{0,\TT}\leq 2M^k$ and $|{\vec n}_\TT|\leq 2M^k$. This part of the proof 
is the most challenging one. Given $\rho_\TT\in {\mathcal S}^k$ and 
$n_\TT=F_\mu^k(\rho_\TT)$, we will first establish the positivity of $n_{\pm,\TT}$ 
(Step 2), then the upper bounds for $n_{\pm,\TT}$ (Step 3), and finally the 
$L^{\infty}$ bound for ${\vec n}_{\perp,\TT}$ (Step 4). 

\medskip
{\em Step 1: Existence and uniqueness of solutions to \eqref{ex.n0}-\eqref{ex.bc2}.}
The linear system of equations \eqref{ex.n0}-\eqref{ex.bc2} is a square system of 
size $4\theta$. The existence of a solution is equivalent to the uniqueness of a 
solution and to the invertibility of the corresponding matrix. Therefore, we just 
have to prove that if the right-hand side to the system is zero then the solution  
is zero. Thus, we may work with the original linear system
assuming homogeneous Dirichlet boundary conditions and setting
$n_{0,K}^{k-1}=\rho_{0,K}=0$ and $\vec n_K^{k-1}=\vec\rho_K=0$, in order 
to set the right-hand side to zero.

We multiply the corresponding equation
\eqref{ex.n0} by $\frac14 n_{0,K}$ and \eqref{ex.nl}
by $n_{\ell,K}$, sum these four equations, and sum over all control volumes $K\in\TT$:
\begin{align*}
  0 &= \left(1+\mu\right)
	\sum_{K\in\TT}\frac{\m(K)}{4\triangle t}n_{0,K}^2
	+ \left(1+\mu\right)
	\sum_{K\in\TT}\frac{\m(K)}{\triangle t}|\vec n_{K}|^2
	+ \frac14\sum_{K\in\TT}\sum_{\sigma\in\EE_K}j_{0,K,\sigma}n_{0,K} \\
	&\phantom{xx}{}+ \sum_{K\in\TT}\sum_{\sigma\in\EE_K}\vec{j}_{K,\sigma}\cdot\vec n_K
	- 2\gamma\sum_{K\in\TT}\m(K)(\vec n_K\times\vec m)\cdot \vec n_K
	+ \frac{\m(K)}{\tau}\sum_{K\in\TT}|\vec n_K|^2 \\
	&= T_1 + \cdots + T_6.
\end{align*}
Note that $T_5=0$ and $T_1$, $T_2$, and $T_6$ are nonnegative.
Thus, it remains to estimate the terms $T_3$ and $T_4$. 
By discrete integration by parts
(note that the problem is homogeneous) and the definitions 
\eqref{sch.j0}-\eqref{sch.jl} of $j_{\ell,K,\sigma}$ (omitting the superindex $k$),
\begin{align*}
  T_3 &= -\frac{D}{4\eta^2}\sum_{\sigma\in\EE}(J_{0,K,\sigma}
	- 2p\vec J_{K,\sigma}\cdot\vec m)\DD n_{0,K,\sigma} =: T_{31} + T_{32}, \\
	T_4 &=- \frac{D}{\eta^2}\sum_{\sigma\in\EE}\left(\eta\vec J_{K,\sigma}
	+ (1-\eta)(\vec J_{K,\sigma}\cdot\vec m)\vec m - \frac{p}{2}J_{0,K,\sigma}\vec m
	\right)\cdot\DD\vec n_{K,\sigma} =: T_{41} + T_{42} + T_{43}.
\end{align*}
With formulation \eqref{sch.coth}, definition \eqref{Bs} of $B^s$, 
and the discrete chain rule $(n_{\ell,K}+n_{\ell,K,\sigma})$ 
$\times\DD n_{\ell,K,\sigma}=\DD(n_{\ell}^2)_{K,\sigma}$, we have
\begin{align*}
  T_{31} &= \frac{D}{8\eta^2}\sum_{\sigma\in\EE}\tau_\sigma\left(2B^s(\DD V_{K,\sigma})
	(\DD n_{0,K,\sigma})^2
	+ \DD(n_{0}^2)_{K,\sigma}\DD V_{K,\sigma}\right), \\
	T_{32} &= -\frac{pD}{4\eta^2}\sum_{\sigma\in\EE}\tau_\sigma
	\big(2B^s(\DD V_{K,\sigma})\DD \vec{n}_{K,\sigma}\cdot\vec{m} 
	+ (\vec{n}_K+\vec{n}_{K,\sigma})\cdot\vec{m}
	\DD V_{K,\sigma}\big)\DD n_{0,K,\sigma}, \\
	T_{41} &= \frac{D}{2\eta}\sum_{\sigma\in\EE}\tau_\sigma\left(2B^s(\DD V_{K,\sigma})
	|\DD\vec{n}_{K,\sigma}|^2
	+ \DD(|\vec{n}|^2)_{K,\sigma}\DD V_{K,\sigma}\right), \\
	T_{42}	&= \frac{(1-\eta)D}{2\eta^2}\sum_{\sigma\in\EE}\tau_\sigma
	\big(2B^s(\DD V_{K,\sigma})
	(\DD\vec{n}_{K,\sigma}\cdot\vec{m})^2 
	+ \DD((\vec{n}\cdot\vec{m})^2)_{K,\sigma}\DD V_{K,\sigma}\big), \\
	T_{43}	&= -\frac{pD}{4\eta^2}\sum_{\sigma\in\EE}\tau_\sigma
	\big(2B^s(\DD V_{K,\sigma})\DD n_{0,K,\sigma}
	+ (n_{0,K}+n_{0,K,\sigma})\DD V_{K,\sigma}\big)\DD\vec{n}_{K,\sigma}\cdot\vec{m}.
\end{align*}

We collect all terms from $T_3+T_4$ involving the function $B^s$:
\begin{align*}
  I_1 &:= \frac{D}{\eta^2}\sum_{\sigma\in\EE}\tau_\sigma B^s(\DD V_{K,\sigma})
  \bigg(\frac14 (\DD n_{0,K,\sigma})^2
	+ \eta|\DD\vec{n}_{K,\sigma}|^2 + (1-\eta)(\DD\vec{n}_{K,\sigma}\cdot\vec{m})^2 
	\\
	&\phantom{xx}{}- p\DD n_{0,K,\sigma}\DD\vec{n}_{K,\sigma}\cdot\vec{m}\bigg) \\
	&= \frac{D}{\eta^2}\sum_{\sigma\in\EE}\tau_\sigma B^s(\DD V_{K,\sigma})\bigg[
	\begin{pmatrix} \DD n_{0,K,\sigma} \\ \DD\vec{n}_{K,\sigma}\cdot\vec{m} 
	\end{pmatrix}^\top \begin{pmatrix} 1/4 & -p/2 \\ -p/2 & 1-\eta^2/2 \end{pmatrix}
	\begin{pmatrix} \DD n_{0,K,\sigma} \\ \DD\vec{n}_{K,\sigma}\cdot\vec{m} 
	\end{pmatrix} \\
	&\phantom{xx}{}+ \frac{\eta^2}{2}|\DD \vec{n}_{K,\sigma}|^2
	+ \eta\left(1-\frac{\eta}{2}\right)\big(|\DD \vec{n}_{K,\sigma}|^2
	- |\DD\vec{n}_{K,\sigma}\cdot\vec{m}|^2\big)\bigg]. 
\end{align*}
The eigenvalues of the $2\times 2$ matrix appearing in $I_1$ are 
$$
  \lambda_\pm = \frac18(5-2\eta^2) \pm \frac18\sqrt{(5-2\eta^2)^2-8\eta^2} 
	\ge \frac14 > 0.
$$
Then, using the inequalities $B^s(z)\ge 1$ for all $z\in\R$ and 
$|\DD\vec{n}_{K,\sigma}|^2 \ge |\DD\vec{n}_{K,\sigma}\cdot\vec{m}|^2$
(since $|\vec{m}|^2=1$), it follows that
\begin{align*}
  I_1 &\ge \frac{D}{\eta^2}\sum_{\sigma\in\EE}\tau_\sigma\left(
	\lambda_-\big((\DD n_{0,K,\sigma})^2 
	+ (\DD \vec{n}_{K,\sigma}\cdot\vec{m})^2\big) 
	+ \frac{\eta^2}{2}|\DD \vec{n}_{K,\sigma}|^2\right) \geq 0.
%	&\ge \frac{D}{\eta^2}\min\left(\lambda_-,\frac{\eta^2}{2}\right)
%	(|n_{0,\MM}|_{1,\MM}^2 + |\vec n_\MM|_{1,\MM}^2),
\end{align*}
%where 
%$$
%  \lambda_\pm = \frac18(5-2\eta^2) \pm \frac18\sqrt{(5-2\eta^2)^2-8\eta^2} 
%	\ge \frac14 > 0
%$$
%are the eigenvalues of the $2\times 2$ matrix appearing in $I_1$.
Next, we collect in $I_2$ the remaining terms from $T_3+T_4$ involving the discrete gradient 
$\DD V_{K,\sigma}$. Taking into account that 
$$
  (\vec{n}_K+\vec{n}_{K,\sigma})\cdot\vec{m}\DD n_{0,K,\sigma}
	+ (n_{0,K}+n_{0,K,\sigma})\DD \vec{n}_{K,\sigma}\cdot\vec{m}
	= 2\DD((\vec{n}\cdot\vec{m})n_0)_{K,\sigma},
$$
integrating by parts, and employing the discrete Poisson equation \eqref{ex.v}, 
we infer that
\begin{align*}
  I_2 &:= \frac{D}{2\eta^2}\sum_{\sigma\in\EE}\tau_\sigma \DD V_{K,\sigma}
	\bigg(\frac14\DD(n_0^2)_{K,\sigma} + \eta\DD(|\vec{n}|^2)_{K,\sigma}
	+ (1-\eta)\DD((\vec{n}\cdot\vec{m})^2)_{K,\sigma} \\
	&\phantom{xx}{}- p\DD((\vec{n}\cdot\vec{m})n_0)_{K,\sigma}\bigg) \\
	&= \frac{D}{2\eta^2\lambda_D^2}\sum_{K\in\TT}\m(K)(\rho_{0,K}-C_K)
	\bigg(\frac14 n_{0,K}^2 + \eta|\vec{n}_K|^2 + (1-\eta)(\vec{n}_K\cdot\vec{m})^2 \\
  &\phantom{xx}{}- p(\vec{n}_K\cdot\vec{m})n_{0,K}\bigg).
\end{align*}
The sum of the terms in the brackets is nonnegative since
\begin{align*}
  \frac14 n_{0,K}^2 &+ \eta|\vec{n}_K|^2 + (1-\eta)(\vec{n}_K\cdot\vec{m})^2 
  - p(\vec{n}_K\cdot\vec{m})n_{0,K} \\
	&= \left(\frac12 n_{0,K}-p\vec{n}_K\cdot\vec{m}\right)^2 
	+ \eta^2(\vec{n}_K\cdot\vec{m})^2 
	+ \eta\big(|\vec{n}_K|^2 - (\vec{n}_K\cdot\vec{m})^2\big) \ge 0.
\end{align*}
Therefore, the term involving $\rho_{0,K}\ge 0$ can be omitted, giving
\begin{align*}
  I_2 &\ge -\frac{D}{2\eta^2\lambda_D^2}\sum_{K\in\TT}\m(K)C_K
	\bigg(\frac14 n_{0,K}^2 + \eta|\vec{n}_K|^2 + (1-\eta)(\vec{n}_K\cdot\vec{m})^2  
	- p(\vec{n}_K\cdot\vec{m})n_{0,K}\bigg) \\
	&\ge -\frac{D}{2\eta^2\lambda_D^2}\|C\|_\infty
	\sum_{K\in\TT}\m(K)\left(\frac12 n_{0,K}^2 + 2|\vec{n}_K|^2\right)
	= -\frac{D}{\eta^2\lambda_D^2}\|C\|_\infty
	\left(\frac14 \|n_{0,\TT}\|_{2}^2 + \|\vec n_\TT\|_{2}^2\right).
\end{align*}

This shows finally that
$$
  T_3 + T_4 \ge 
	- \frac{D\|C\|_\infty}{\eta^2\lambda_D^2}\left(\frac14\|n_{0,\TT}\|_{2}^2 
	+ \|\vec n_\TT\|_{2}^2\right),
$$
and summarizing the estimates for $T_1,\ldots,T_6$, we conclude that
$$
 \left(\frac{1+\mu}{\triangle t}
	-\frac{D\|C\|_\infty}{\eta^2\lambda_D^2}\right)\left(\frac14\|n_{0,\TT}\|_{2}^2 
	+ \|\vec n_\TT\|_{2}^2\right) \le 0.
$$
Hence, choosing $\mu$ as in \eqref{ex.mu}, the first bracket is positive, 
showing that $n_{0,\TT}=0$ and
$\vec n_\TT=0$, which proves the invertibility of the linear system of 
equations \eqref{ex.n0}-\eqref{ex.bc2}. 
The second step involved in the definition of $F_\mu^k$,  
$(V_\TT,\rho_\TT) \to n_\TT$, is a well-defined mapping. 
Moreover, the matrix and the right-hand side of the linear system of equations 
are continuous with respect to $(V_\TT, \rho_\TT)$ so that the mapping is continuous.
\medskip

%%%%%%%%%%%%%%%%%%%%%%%%%%%%%%%%%%%%%%%%%%%
{\em Step 2: Positivity of $n_{\pm,\TT}$.}
We will prove that $n_{\pm,K}\ge 0$ for all $K\in\TT$.
Multiplying \eqref{ex.n0}-\eqref{ex.nl} by $\vec m$ and adding or subtracting it
from \eqref{ex.n0}, multiplied by $\frac12$, we find that
\begin{align}\label{ex.npm} 
  \frac{\m(K)}{\triangle t} & (n_{\pm,K} - n_{\pm,K}^{k-1})
	+ \mu\frac{\m(K)}{\triangle t}(n_{\pm,K}-\rho_{\pm,K})
	+ D(1\pm p)\sum_{\sigma\in\EE_K}J_{\pm,K,\sigma} \\ 
	&= \mp\frac{\m(K)}{2\tau}(n_{+,K}-n_{-,K}), \nonumber
\end{align}
where $\rho_{\pm,K}=\frac12 \rho_{0,K}\pm\vec\rho_K\cdot\vec m_K$ and 
$J_{\pm,K,\sigma}=\frac12 J_{0,K,\sigma} \pm {\vec J}_{K,\sigma}\cdot {\vec m}$, i.e.
\begin{equation}\label{def.Jpm}
  J_{\pm,K,\sigma}=\tau_{\sigma} \big( B(\DD V_{K,\sigma}) n_{\pm,K}
	- B(-\DD V_{K,\sigma}) n_{\pm,K,\sigma}\big).
\end{equation}
Then, multiplying \eqref{ex.npm} by $n_{\pm,K}^-=\min\{0,n_{\pm,K}\}$,
summing over all control volumes $K\in\TT$, and adding both equations, 
it follows that
\begin{align}
  0 &= \sum_\pm\sum_{K\in\TT}\frac{\m(K)}{\triangle t}
	(n_{\pm,K}-n_{\pm,K}^{k-1})n_{\pm,K}^-
	+ \mu\sum_\pm\sum_{K\in\TT}
	\frac{\m(K)}{\triangle t}(n_{\pm,K}-\rho_{\pm,K})n_{\pm,K}^- \nonumber \\
	&\phantom{xx}{}+ D\sum_\pm(1\pm p)\sum_{K\in\TT}\sum_{\sigma\in\EE_K}
	J_{\pm,K,\sigma} n_{\pm,K}^- 
	+\sum_{K\in\TT}\frac{\m(K)}{2\tau}(n_{+,K}-n_{-,K})
	(n_{+,K}^--n_{-,K}^-)  \label{ex.aux} \\
	&=: T_7 + T_8 + T_9 + T_{10}, \nonumber
\end{align}
Since $n_{\pm,K} n_{\pm,K}^- = (n_{\pm,K}^-)^2$, 
$n_{\pm,K}^{k-1}\ge 0$, and $\rho_{\pm,K}\ge 0$, 
the first two terms in \eqref{ex.aux} can be estimated as
$$
  T_7 \ge \sum_\pm\sum_{K\in\TT}\frac{\m(K)}{\triangle t}(n_{\pm,K}^-)^2, \quad
	T_8 \ge \mu\sum_\pm\sum_{K\in\TT}\frac{\m(K)}{\triangle t}
	(n_{\pm,K}^-)^2.
$$ 
The monotonicity of the mapping $z\mapsto z^-$ shows that $T_{10}$ is nonnegative.
By the discrete integration-by-parts formula \eqref{ibp}, the third term 
in \eqref{ex.aux} becomes
$$
  T_9 = -D\sum_\pm(1\pm p)\sum_{\sigma\in\EE}J_{\pm,K,\sigma}\DD (n_{\pm}^-)_{K,\sigma}.
$$
The sum over the boundary edges vanishes since $n_{\pm,K,\sigma}^-=0$ for all 
$\sigma\in \EE^D_{\rm ext}$. We claim that 
\begin{equation}\label{ex.claim}
  -J_{\pm,K,\sigma}\DD (n_{\pm}^-)_{K,\sigma}
	\ge \frac{1}{2}\tau_\sigma 
	\DD V_{K,\sigma} \big((n_{\pm,K,\sigma}^-)^2
	- (n_{\pm,K}^-)^2\big),\quad  \mbox{for } K\in\TT\mbox{ and } \sigma \in \EE_K,
\end{equation}
such that
\begin{equation}\label{ex.t9}
  T_9 \ge \frac{D}{2}\sum_\pm(1\pm p)\sum_{\sigma\in\EE}\tau_\sigma
	\DD V_{K,\sigma}\DD \big((n_{\pm}^-)^2\big)_{K,\sigma}.
\end{equation}

To prove \eqref{ex.claim}, we distinguish the cases $\DD V_{K,\sigma}\ge 0$ and 
$\DD V_{K,\sigma}\le 0$. If $\DD V_{K,\sigma}\ge 0$, we apply a formulation 
similar to \eqref{1.refJ.2}, leading to 
$$
 -J_{\pm,K,\sigma}\DD (n_{\pm}^-)_{K,\sigma}
	= \tau_\sigma \big(\DD V_{K,\sigma} n_{\pm,K,\sigma}
	+ B(\DD V_{K,\sigma})\DD n_{\pm,K,\sigma}\big)\DD (n_{\pm}^-)_{K,\sigma}.
$$
Then, using the nonnegativity of the function $B$, the monotonicity of the mapping 
$z\mapsto z^-$, and the inequality $z(z^--y^-)\geq \frac{1}{2}((z^-)^2-(y^-)^2)$, 
we obtain \eqref{ex.claim}. If $\DD V_{K,\sigma}\le 0$,
we employ formulation \eqref{1.refJ.1}, so that 
$$
-J_{\pm,K,\sigma}\DD (n^-_\pm)_{K,\sigma}
	= \tau_\sigma \big(\DD V_{K,\sigma} n_{\pm,K}
	+ B(-\DD V_{K,\sigma})\DD n_{\pm,K,\sigma}\big)\DD (n_{\pm}^-)_{K,\sigma}
$$
and similar arguments lead to \eqref{ex.claim}.

Applying discrete integration by parts to the right-hand side of \eqref{ex.t9} 
(the boundary term vanishes since the boundary data is nonnegative) and 
employing the discrete Poisson equation \eqref{ex.v}, we find that
\begin{align*}
 T_9&\geq  -\frac{D}{2}\sum_{K\in\TT}  \sum_{\sigma\in\EE_K}\tau_\sigma \DD V_{K,\sigma}
	\big((1+p)(n_{+,K}^-)^2 + (1-p)(n_{-,K}^-)^2\big) \\
	&= \frac{D}{2\lambda_D^2}\sum_{K\in\TT}\m(K)(\rho_{0,K}-C_K)
	\big((1+p)(n_{+,K}^-)^2 + (1-p)(n_{-,K}^-)^2\big) \\
	&\ge -\frac{D}{2\lambda_D^2}\|C\|_{\infty}\sum_{K\in\TT}\m(K)
	\big((1+p)(n_{+,K}^-)^2 + (1-p)(n_{-,K}^-)^2\big).
\end{align*}

Summarizing the above estimates, we conclude from \eqref{ex.aux} that
$$
  \left( \frac{1+\mu}{\triangle t}-\frac{D(1+p)}{2\lambda_D^2}\|C\|_{\infty}\right)
  \sum_{K\in\TT}	\m(K)\big((n_{+,K}^-)^2 + (n_{-,K}^-)^2\big)\leq 0.
$$
By the choice of $\mu$ in \eqref{ex.mu}, we deduce that
$$
  \sum_{K\in\TT}\m(K)\big((n_{+,K}^-)^2 + (n_{-,K}^-)^2\big) \le 0,
$$
which implies that $n_{\pm,K}^-=0$ and hence $n_{\pm,K}\ge 0$ for all $K\in\TT$.
\medskip

{\em Step 3: Upper bounds for $n_{\pm,\TT}$.}
The goal is to show that $n_{\pm,K}\le M^k$ for all $K\in\TT$, where 
$M^k$ is defined in Theorem \ref{thm.ex}.
We multiply \eqref{ex.npm} by $(n_{\pm,K}-M^k)^+=\max\{0,n_{\pm,K}-M^k\}$, sum
over all $K\in\TT$, and add both equations:
\begin{align}
  0 &= \sum_\pm\sum_{K\in\TT}\frac{\m(K)}{\triangle t}
	\big((n_{\pm,K}-M^k)-(n_{\pm,K}^{k-1}-M^{k-1})\big)(n_{\pm,K}-M^k)^+  \nonumber \\
	&\phantom{xx}{}+ \mu\sum_\pm\sum_{K\in\TT}
	\frac{\m(K)}{\triangle t}\big((n_{\pm,K}-M^k)-(\rho_{\pm,K}-M^k)\big)
	(n_{\pm,K}-M^k)^+ \nonumber \\
	&\phantom{xx}{}+\sum_\pm \sum_{K\in\TT}\frac{\m(K)}{\triangle t}(M^k-M^{k-1})
	(n_{\pm,K}-M^k)^+ \label{ex.aux3} \\
	&\phantom{xx}{}
	+ D\sum_\pm(1\pm p)\sum_{K\in\TT}\sum_{\sigma\in\EE_K}J_{\pm,K,\sigma}
	(n_{\pm,K}-M^k)^+ \nonumber \\
	&\phantom{xx}{}+\frac{1}{2\tau}\sum_\pm\sum_{K\in\TT}\m(K)
	\big((n_{+,K}-M^k)-(n_{-,K}-M^k)\big)(\pm(n_{\pm,K}-M^k)^+) \nonumber \\
	&=: T_{11} + T_{12} + T_{13} + T_{14} + T_{15}. \nonumber 
\end{align}
Using the inequality $(z-y)z^+ \ge \frac12((z^+)^2-(y^+)^2)$,
the first two terms are estimated by
\begin{align*}
  T_{11} &\ge 
  \frac{1}{2\triangle t}\sum_\pm\sum_{K\in\TT}
	\m(K)\big((n_{\pm,K}-M^k)^+\big)^2, \\
	T_{12} &\ge \frac{\mu}{2\triangle}\sum_\pm\sum_{K\in\TT}
	\m(K)\big((n_{\pm,K}-M^k)^+\big)^2,
\end{align*}
since $n_{\pm,K}^{k-1}\le M^{k-1}$ and $\rho_{\pm,K}\le M^k$ by assumption. 
By definition of $M^k$, the third term $T_{13}$ becomes
$$
  T_{13} = \alpha M^k\sum_\pm\sum_{K\in\TT}\m(K)(n_{\pm,K}-M^k)^+,
$$
and the last term $T_{15}$ is nonnegative.

It remains to estimate $T_{14}$.
By discrete integration by parts (the boundary term vanishes in view of
$n_{\pm,\sigma}^D\le M^k$ for $\sigma\in\EE_{\rm ext}^D$), we find that
$$
  T_{14} = -D\sum_\pm(1\pm p)\sum_{\sigma\in\EE}J_{\pm,K,\sigma}
	\DD\big((n_\pm-M^k)^+\big)_{K,\sigma}.
$$
Similarly as in Step 2, we claim that the following estimate holds: 
$$
  T_{14} \ge D\sum_\pm(1\pm p)\sum_{\sigma\in\EE}\tau_\sigma\DD V_{K,\sigma}
	\DD \left(\frac12((n_\pm-M^k)^+)^2 + M^k(n_\pm-M^k)^+\right)_{K,\sigma}.
$$
Indeed, let first $\DD V_{K,\sigma}\ge 0$. Using the inequalities
$\DD(n_{\pm}-M^k)_{K,\sigma}\DD((n_\pm-M^k)^+)_{K,\sigma}\ge 0$ and
$(n_{\pm,K,\sigma}-M^k)\DD((n_\pm-M^k)^+)_{K,\sigma}\ge 
\frac12\DD(((n_\pm-M^k)^+)^2)_{K,\sigma}\ge 0$, it follows from \eqref{1.refJ.2} that
\begin{align*}
  -J_{\pm,K,\sigma} 
	&\ge \tau_\sigma\DD V_{K,\sigma}\big((n_{\pm,K,\sigma}-M^k)+M^k\big)
	\DD((n_\pm-M)^+)_{K,\sigma} \\
	&\ge \tau_\sigma\DD V_{K,\sigma}\left(\frac12
	\DD(((n_\pm-M^k)^+)^2)_{K,\sigma}+M^k \DD((n_\pm-M^k)^+)_{K,\sigma}\right).
\end{align*}
The proof for $\DD V_{K,\sigma}\le 0$ is similar, employing formulation
\eqref{1.refJ.1}. Then, integrating by parts and employing the Poisson equation
and $\rho_{0,K}\ge 0$,
\begin{align*}
  T_{14} &\ge \frac{D}{\lambda^2}\sum_\pm(1\pm p)\sum_{K\in\TT}\m(K)(\rho_{0,K}-C_K)
	\left(\frac12((n_{\pm,K}-M^k)^+)^2 + M^k(n_{\pm,K}-M^k)^+\right) \\
  &\ge -\frac{D}{\lambda_D^2}\|C\|_\infty(1+p)\sum_\pm\sum_{K\in\TT}\m(K)
	\left(\frac12((n_{\pm,K}-M^k)^+)^2 + M^k(n_{\pm,K}-M^k)^+\right).
\end{align*}

Summarizing the above estimates, we infer from \eqref{ex.aux3} that
\begin{align*}
  & \left(\frac{1+\mu}{2\triangle t}
	- \frac{D}{\lambda_D^2}\|C\|_\infty(1+p)\right)
	\sum_\pm\sum_{K\in\TT}\m(K)\big((n_{\pm,K}-M^k)^+\big)^2 \\
	&\phantom{xx}{}+ \left(\alpha - \frac{D}{\lambda_D^2}\|C\|_\infty(1+p)\right) 
	M^k\sum_\pm\sum_{K\in\TT}\m(K)(n_{\pm,K}-M^k)^+
	\le 0.
\end{align*}
Then, choosing $\mu$ as in \eqref{ex.mu} and taking into account the
definition of $\alpha$, we infer that $n_{\pm,K}\le M^k$ for $K\in\TT$. 

%%%%%%%%%%%%%%%%%%%%
\medskip

{\em Step 4: $L^\infty$ bound for $\vec n_{\perp,\TT}^k$.}
We prove a uniform $L^{2q}$ bound for $\vec n_{\perp,K}
=\vec n_K-(\vec n_K\cdot\vec m)\vec m$. For this,
we multiply the vector version of \eqref{ex.nl} (omitting the superindex $k$)
by $\vec{m}$ twice, and taking the difference of the equations for $\vec{n}_K$ and 
$(\vec{n}_K\cdot\vec{m})\vec{m}$, we obtain
\begin{align*}
  \frac{\m(K)}{\triangle t} & \big(\vec{n}_{\perp,K} - \vec{n}_{\perp,K}^{k-1}
	+ \mu(\vec n_{\perp,K}-\vec{\rho}_{\perp,K}) \big)
	+ \frac{D}{\eta}\sum_{\sigma\in\EE_K}\vec{J}_{\perp,K,\sigma} 
	- 2\gamma\m(K)(\vec{n}_{\perp,K}\times\vec{m}) \\
	&= -\frac{\m(K)}{\tau}\vec{n}_{\perp,K},
\end{align*}
where $\vec{\rho}_{\perp,K}=\vec\rho_K-(\vec\rho_K\cdot\vec m)\vec m$,
and $\vec J_{\perp,K,\sigma}$ is given by
\begin{align}
  \vec J_{\perp,K,\sigma} &= \tau_\sigma\big(-\DD V_{K,\sigma}\vec n_{\perp,K}
	- B(-\DD V_{K,\sigma})\DD \vec n_{\perp,K,\sigma}\big) \label{ex.Jperp.1}\\
	&= \tau_\sigma\big(-\DD V_{K,\sigma}\vec n_{\perp,K,\sigma}
	- B(\DD V_{K,\sigma})\DD \vec n_{\perp,K,\sigma}\big) . \label{ex.Jperp.2}
\end{align}
Then, multiplying this equation
by $|\vec{n}_{\perp,K}|^{2(q-1)}\vec{n}_{\perp,K}$ (where $q\in\N$)
and summing over $K\in\TT$, we arrive at $T_{16}+T_{17}+T_{18} + T_{19}=0$ with 
\begin{align*}
  T_{16} &= \frac{1}{\triangle t}\sum_{K\in\TT} \m(K)
	(\vec{n}_{\perp,K}-\vec{n}_{\perp,K}^{k-1})
	\cdot \vec{n}_{\perp,K}|\vec{n}_{\perp,K}|^{2(q-1)},\\
  T_{17}&=\frac{\mu}{\triangle t}
	\sum_{K\in\TT}\m(K)
	(\vec{n}_{\perp,K}-\vec{\rho}_{\perp,K})
	\cdot \vec{n}_{\perp,K}|\vec{n}_{\perp,K}|^{2(q-1)},\\
  T_{18}&= \frac{D}{\eta}
	\sum_{K\in\TT}\sum_{\sigma\in\EE_K}\vec{J}_{\perp,K,\sigma}\cdot
	\vec{n}_{\perp,K}|\vec{n}_{\perp,K}|^{2(q-1)},\\
  T_{19}&=\frac{1}{\tau}\sum_{K\in\TT}\m(K)|\vec{n}_{\perp,K}|^{2q}. 
\end{align*}
The elementary inequality $|\vec a|^{2(q-1)}\vec a\cdot(\vec a-\vec b)
\ge (|\vec a|^{2q}-|\vec b|^{2q})/(2q)$ for $\vec a$, $\vec b\in\R^3$
shows that
\begin{align*}
  T_{16} &\ge \frac{1}{2q\triangle t}
	\left(\sum_{K\in\TT}\m(K)|\vec{n}_{\perp,K}|^{2q}
	- \sum_{K\in\TT}\m(K)|\vec{n}_{\perp,K}^{k-1}|^{2q}\right), \\
  T_{17} &\ge \frac{\mu}{2q\triangle t}
	\left(\sum_{K\in\TT}\m(K)|\vec{n}_{\perp,K}|^{2q}
	- \sum_{K\in\TT}\m(K)|\vec{\rho}_{\perp,K}|^{2q}\right).
\end{align*}
By discrete integration by parts (observe that $\vec n_{\perp,\sigma}=0$
for $\sigma\in\EE_{\rm ext}^D$),
\begin{align*}
  T_{18} 
	&= -\frac{D}{\eta}\sum_{\sigma\in\EE}\vec{J}_{\perp,K,\sigma}\cdot
	\big(\vec{n}_{\perp,K,\sigma}|\vec{n}_{\perp,K,\sigma}|^{2(q-1)}
	- \vec{n}_{\perp,K}|\vec{n}_{\perp,K}|^{2(q-1)}\big).
\end{align*}
Again, we distinguish the cases $\DD V_{K,\sigma}\ge 0$
and $\DD V_{K,\sigma}<0$ for given $K\in\TT$ and $\sigma\in\EE_K$.
First, let $\DD V_{K,\sigma}\ge 0$ and use formulation \eqref{ex.Jperp.2} 
of the numerical flux. This gives
\begin{align*}
  -\vec{J}_{\perp,K,\sigma} & \cdot
	\big(\vec{n}_{\perp,K,\sigma}|\vec{n}_{\perp,K,\sigma}|^{2(q-1)}
	- \vec{n}_{\perp,K}|\vec{n}_{\perp,K}|^{2(q-1)}\big) \\
	&= \tau_\sigma\DD V_{K,\sigma}\vec{n}_{\perp,K,\sigma}\cdot
	\big(\vec{n}_{\perp,K,\sigma}|\vec{n}_{\perp,K,\sigma}|^{2(q-1)}
	- \vec{n}_{\perp,K}|\vec{n}_{\perp,K}|^{2(q-1)}\big) \\
	&\phantom{xx}{}+ B(\DD V_{K,\sigma})(\vec{n}_{\perp,K,\sigma}-\vec{n}_{\perp,K})
	\cdot\big(\vec{n}_{\perp,K,\sigma}|\vec{n}_{\perp,K,\sigma}|^{2(q-1)}
	- \vec{n}_{\perp,K}|\vec{n}_{\perp,K}|^{2(q-1)}\big) 
\end{align*}
Because of 
\begin{align*}
  \vec{a}\cdot(\vec{a}|\vec{a}|^{2(q-1)}-\vec{b}|\vec{b}|^{2(q-1)})
	&= |\vec{a}|^{2q} - \vec{a}\cdot\vec{b}|\vec{b}|^{2(q-1)} 
	\ge |\vec{a}|^{2q} - \frac{1}{2q}|\vec{a}|^{2q} - \left(1-\frac{1}{2q}\right)
	|\vec{b}|^{2q} \\
	&\geq \left(1-\frac{1}{2q}\right)(|\vec{a}|^{2q}-|\vec{b}|^{2q}) 
	\quad\mbox{for all }\vec{a},\vec{b}\in\R^3,
\end{align*}
applied to $\vec{a}=\vec{n}_{\perp,K,\sigma}$ and $\vec{b}=\vec{n}_{\perp,K}$, 
and the monotonicity of the mapping $\vec{a}\mapsto \vec{a}|\vec{a}|^{2(q-1)}$,
we find that
\begin{align*}
  -\vec{J}_{\perp,K,\sigma} & \cdot
	\big(\vec{n}_{\perp,K,\sigma}|\vec{n}_{\perp,K,\sigma}|^{2(q-1)}
	- \vec{n}_{\perp,K}|\vec{n}_{\perp,K}|^{2(q-1)}\big) \\
  &\ge \tau_\sigma\left(1-\frac{1}{2q}\right)\DD V_{K,\sigma}
	\big(|\vec{n}_{\perp,K,\sigma}|^{2q}-|\vec{n}_{\perp,K}|^{2q}\big).
\end{align*}
This result still holds if $\DD V_{K,\sigma}<0$, thanks to formulation 
\eqref{ex.Jperp.1}. Therefore
$$
  T_{18} \geq \frac{D}{\eta}\left(1-\frac{1}{2q}\right)
	\sum_{\sigma\in\EE} \tau_\sigma\DD V_{K,\sigma}\DD(|\vec{n}_\perp|^{2q})_{K,\sigma}.
$$
Using discrete integration by parts and the Poisson equation \eqref{sch.v} leads to
\begin{align*}
  T_{18} &\ge \frac{D}{\eta\lambda_D^2}\left(1-\frac{1}{2q}\right)
	\sum_{K\in\TT}\m(K)(\rho_{0,K}-C_K)|\vec{n}_{\perp,K}|^{2q} 
	\ge -\frac{D}{\eta\lambda_D^2}\|C\|_{\infty}
	\sum_{K\in\TT}\m(K)|\vec{n}_{\perp,K}|^{2q}.
\end{align*}

Summarizing the above estimates, we obtain
\begin{align*}
  \left(1 + \mu + 2q\triangle t\Big(\frac{1}{\tau}
	-\frac{D\|C\|_{\infty}}{\eta\lambda_D^2}\Big)\right)
	\sum_{K\in\TT}\m(K)|\vec{n}_{\perp,K}|^{2q}
	&\le \sum_{K\in\TT}\m(K)|\vec{n}_{\perp,K}^{k-1}|^{2q} \\
	&\phantom{xx}{}{}+ \mu\sum_{K\in\TT}\m(K)
	|\vec{\rho}_{\perp,K}|^{2q} .
\end{align*}
Condition \eqref{hypo5} on $\tau$, the induction hypothesis 
$\|\vec n_{\perp,\TT}^{k-1}\|_\infty\le M^{k-1}
\le M^k$, and the fact that $\rho\in {\mathcal S}^k$ (see \eqref{def.Sk} for the
definition of ${\mathcal S}^k$), such that 
$\|\vec\rho_{\perp,\TT}\|_\infty\le M^k$, imply that 
$$
  \|\vec n_{\perp,\TT}\|_{2q}\leq \mbox{meas}(\Omega)^{1/(2q)} M^k
	\quad\mbox{for } q\geq 1.
$$
Passing to the limit $q\to +\infty$, we deduce that 
$\|n_{\perp,\TT}\|_\infty \le M^k$.

%%%%%%%%%%%%%%%%%%%%
\medskip

{\em Conclusion.} In Step 1, we have proved that the mapping $F_\mu^k$ is 
well-defined and continuous. In Steps 2-4, we have proved that $F_\mu^k$ 
preserves the bounded set ${\mathcal S}^k$.
Thus, the fixed-point theorem of Brouwer shows the existence
of a fixed point to $F_\mu^k$, belonging to ${\mathcal S}^k$. Let us denote 
this fixed point by $n_\TT^k=(n_{0,\TT}^k, {\vec n}_\TT^k)$. It is a solution 
to scheme \eqref{sch.n0}--\eqref{sch.bcl} at step $k$ and satisfies 
$$
  0\leq n_{\pm,K}^k\leq M^k \quad\mbox{and}\quad 
	|{\vec n}_{\perp,K}^k|\leq M^k,\quad \mbox{for } K\in \TT.
$$

%%%%%%%%%%%%%%%%%%%%%%%%%%%%%%%%%%%%%%%%%%%

\subsection{Uniform bounds for the spin-up and spin-down densities}
In order to conclude the proof of Theorem \ref{thm.ex}, it remains to prove that 
the upper bounds on the spin-up and spin-down densities in fact do not depend on $k$. 
The positivity of these densities is already proved above.  

We assume as induction hypothesis that $n_{\pm,K}^{k-1}\leq M^0$ for all $K\in\TT$ 
(this property is ensured for $k=1$ by the definition of $M^0$). Scheme 
\eqref{sch.n0}-\eqref{sch.nl} implies that 
\begin{equation}\label{ex.npm.k}
  \frac{\m(K)}{\triangle t} (n_{\pm,K}^k - n_{\pm,K}^{k-1})
	+ D(1\pm p)\sum_{\sigma\in\EE_K}J_{\pm,K,\sigma}^k 
	= \mp\frac{\m(K)}{2\tau}(n_{+,K}^k-n_{-,K}^k).
\end{equation}
As in Step 3 above, we multiply \eqref{ex.npm.k} by $(n_{\pm,K}^k-M^0)^+$, sum 
over all $K\in\TT$ and add both equations. This yields $S_1+S_2+S_3=0$, where 
\begin{align*}
  S_1&=\sum_\pm\sum_{K\in\TT}\frac{\m(K)}{\triangle t}
	\big((n_{\pm,K}^k-M^0)-(n_{\pm,K}^{k-1}-M^{0})\big)(n_{\pm,K}^k-M^0)^+, \\
  S_2&=D\sum_\pm (1\pm p)\sum_{K\in\TT}\sum_{\sigma\in\EE_K}J_{\pm,K,\sigma}^k 
	(n_{\pm,K}^k-M^0)^+, \\
  S_3&=\frac{1}{2\tau}\sum_{K\in\TT} \m (K) (n_{+,K}^k-n_{-,K}^k)
	\left((n_{+,K}^k-M^0)^+-(n_{-,K}^k-M^0)^+\right).
\end{align*}
It is clear that $S_3\geq 0$ and, by the induction hypothesis, that 
$$
  S_1\geq \frac{1}{2\triangle t} \sum_\pm\sum_{K\in\TT} \m(K) 
	\left((n_{\pm,K}^k-M^0)^+\right)^2.
$$
The term $S_2$ is the analogue of $T_{14}$. Following the same ideas as in Step 3, 
we obtain
$$
  S_2 \ge \frac{D}{\lambda^2}\sum_\pm(1\pm p)\sum_{K\in\TT}\m(K)(n_{0,K}^k-C_K)
	\left(\frac12((n_{\pm,K}^k-M^0)^+)^2 + M^0(n_{\pm,K}^k-M^0)^+\right).
$$
But, as $n_{0,K}^k=n_{+,K}^k+n_{-,K}^k$, the positivity of $n_{+,K}^k$ and 
$n_{-,K}^k$ and the definition of $M^0$ ensure that 
$n_{0,K}^k-C_K\geq n_{+,K}^k-M^0$ and $n_{0,K}^k-C_K\geq n_{-,K}^k-M^0$, 
leading to $S_2\geq 0$. Therefore, we infer that
$$
  \sum_\pm\sum_{K\in\TT} \m (K) \left((n_{\pm,K}^k-M^0)^+\right)^2\leq 0,
$$
which yields the expected result.

%%%%%%%%%%%%%%%%%%%%%%%%%%%%%%%%%%%%%%%%%%%%%%%%%%%%%%%%%%%%%%%%%%%%%%%%%%%%%%%

\section{Proof of Theorem \ref{thm.diss}}\label{sec.proof2}

Let $(n_{\pm,\TT}^k,V_\TT^k)_{k\geq 0}$ be a solution to \eqref{sch.v}, 
\eqref{ex.npm.k} 
with the corresponding Dirichlet-Neumann boundary conditions.
Since we have to deal with the 
logarithm of the densities $n_{\pm,K}^k$, which may vanish, we introduce a
 regularization of the discrete free energy. 
For $\delta>0$, we set $n_{\pm,K}^{k,\delta}=n_{\pm,K}^k+\delta$ and define 
\begin{align}\label{ex.E}
  E_\delta^k &= \sum_\pm\sum_{K\in\TT}\m(K)\left(n_{\pm,K}^{k,\delta}
	(\log(n_{\pm,K}^{k,\delta})-1) 
	- n_{\pm,K}^k\log\left(\frac{n_{K}^D}{2}+\delta\right)+\frac{n_{K}^D}{2} 
	\right) \\
	&\phantom{xx}{}
	+ \frac{\lambda_D^2}{2}\sum_{\sigma\in\EE}\tau_\sigma(\DD (V^k-V^D)_{K,\sigma})^2. 
	\nonumber 
\end{align}
Therefore, we have  $E_\delta^k-E_\delta^{k-1}=U_1+U_2$, where
\begin{align*}
  U_1 &= \sum_\pm\sum_{K\in\TT}\m(K)\bigg(
	n_{\pm,K}^{k,\delta}(\log n_{\pm,K}^{k,\delta}-1) - n_{\pm,K}^{k-1,\delta}
	(\log n_{\pm,K}^{k-1,\delta}-1) \\
	&\phantom{xx}{}- (n^k_{\pm,K}-n_{\pm,K}^{k-1})
	\log\left(\frac{n_{K}^D}{2}+\delta\right)\bigg), \\
  U_2	& = \frac{\lambda_D^2}{2}\sum_{\sigma\in\EE}\tau_\sigma
	\big((\DD(V^k-V^D)_{K,\sigma})^2 - (\DD(V^{k-1}-V^D)_{K,\sigma})^2\big). 
\end{align*}
The convexity of $x\mapsto x(\log x-1)$ shows that
$x(\log x-1)-y(\log y-1)\le (x-y)\log x$ for all $x$, $y>0$. Hence,
$$
  U_1\leq \sum_\pm\sum_{K\in\TT}\m(K)(n^k_{\pm,K}-n_{\pm,K}^{k-1})\left(
  \log n_{\pm,K}^{k,\delta} -\log\left(\frac{n_{K}^D}{2}+\delta\right)\right).
$$
Using the elementary inequality $\frac12(x^2-y^2)\le (x-y)x$
for all $x$, $y\in\R$, integrating by parts, and employing the
discrete Poisson equation \eqref{sch.v}, it follows that
\begin{align*}
  U_2
	&\le \lambda_D^2\sum_{\sigma\in\EE}\tau_\sigma
	\DD(V^k-V^{k-1})_{K,\sigma}\DD(V^k-V^D)_{K,\sigma} \\
	&= -\lambda_D^2\sum_{K\in\TT}\sum_{\sigma\in\EE_K}\tau_\sigma
	\DD(V^k-V^{k-1})_{K,\sigma}(V_K-V^D_K) \\
	&= \sum_\pm \sum_{K\in\TT}\m(K)
	\big(n_{\pm,K}^k-n_{\pm,K}^{k-1}\big)(V_K^k-V^D_K).
\end{align*}

We summarize the above inequalities and use scheme \eqref{ex.npm.k} to find that 
\begin{align*}
  \frac{1}{\triangle t} & (E_\delta^k-E_\delta^{k-1})
	\leq -\frac{1}{\tau}\sum_{K\in\TT} (n_{+,K}^k-n_{-,K}^k)
	\big(\log n_{+,K}^{k,\delta}-\log n_{-,K}^{k,\delta}\big) \\
  &{}-\sum_\pm D(1\pm p)\sum_{K\in\TT}\sum_{\sigma\in\EE_K}
  J_{\pm,K,\sigma}^k \left(\log n_{\pm,K}^{k,\delta} + V_K^k 
	- \log\left(\frac{n_{K}^D}{2}+\delta\right) - V_K^D\right).
\end{align*}
The first term on the right-hand side is clearly nonpositive. We apply
the discrete inte\-gra\-tion-by-parts formula \eqref{ibp} to the second term. Then,
with the hypothesis on the boundary data (i.e.\ $\log(n^D/2)+V^D$ is 
constant in $\overline\Omega$ such that $\DD V^D_{K,\sigma}
=-\DD (\log n^D)_{K,\sigma}$ for all $K\in\TT$ and $\sigma\in\EE_K$), we infer that
\begin{align*}
  \frac{1}{\triangle t} (E_\delta^k-E_\delta^{k-1})
	&\leq\sum_\pm D(1\pm p)\sum_{\sigma\in\EE}
  J_{\pm,K,\sigma}^k \DD(\log n_{\pm}^{k,\delta}+V^k)_{K,\sigma}\\
  &\phantom{xx}{} + \sum_\pm D(1\pm p)\sum_{\sigma\in\EE}
  J_{\pm,K,\sigma}^k\DD\big(\log n^D-\log (n^D+2\delta)\big)_{K,\sigma}.
\end{align*}
Introducing the numerical fluxes associated to the regularized densities,
$$
  J_{\pm,K,\sigma}^{k,\delta} 
	= \tau_\sigma\big(B(\DD V_{K,\sigma}^k)n_{\pm,K}^{k,\delta}
	- B(-\DD V_{K,\sigma}^k)n_{\pm,K,\sigma}^{k,\delta}\big) 
	= J_{\pm,K,\sigma}^k - \delta\tau_\sigma \DD V_{K,\sigma}^k,
$$
we can write
$$
  \frac{1}{\triangle t} (E_\delta^k-E_\delta^{k-1})\leq U_3+U_4+U_5+U_6,
$$
where
\begin{align*}
  U_3&= \sum_\pm D(1\pm p)\sum_{\sigma\in\EE}
  J_{\pm,K,\sigma}^{k,\delta} \DD\big(\log n_{\pm}^{k,\delta}+V^k\big)_{K,\sigma},\\
  U_4&=\sum_\pm D(1\pm p)\sum_{\sigma\in\EE}
  J_{\pm,K,\sigma}^{k,\delta}\DD(\log n^D-\log (n^D+2\delta))_{K,\sigma},\\
  U_5&=\delta \sum_\pm D(1\pm p)\sum_{\sigma\in\EE}\tau_\sigma \DD V_{K,\sigma}^k
	\DD(\log n_{\pm}^{k,\delta}+V^k)_{K,\sigma},\\
  U_6&=\delta \sum_\pm D(1\pm p)\sum_{\sigma\in\EE}\tau_\sigma \DD V_{K,\sigma}^k
	\DD\big(\log n^D-\log (n^D+2\delta)\big)_{K,\sigma}.
\end{align*}
Now, we employ the following inequalities, which are proved in 
\cite[Appendix A]{BCV14}:
\begin{align*}
  J_{\pm,K,\sigma}^{k,\delta} \DD(\log n_{\pm}^{k,\delta}+V^k)_{K,\sigma}
	&\leq -\tau_{\sigma} \min(n_{\pm,K}^{k,\delta},n_{\pm,K,\sigma}^{k,\delta})
	(\DD(\log n_{\pm}^{k,\delta}+V^k)_{K,\sigma})^2, \\
  \big|J_{\pm,K,\sigma}^{k,\delta}\big| 
	&\leq \tau_{\sigma} \max (n_{\pm,K}^{k,\delta},n_{\pm,K,\sigma}^{k,\delta})
	\big|\DD(\log n_{\pm}^{k,\delta}+V^k)_{K,\sigma}\big|.
\end{align*}
The first inequality yields
$$
  U_3 \leq -\sum_\pm D(1\pm p)\sum_{\sigma\in\EE}\tau_{\sigma} 
	\min(n_{\pm,K}^{k,\delta},n_{\pm,K,\sigma}^{k,\delta})
	\big(\DD(\log n_{\pm}^{k,\delta}+V^k)_{K,\sigma}\big)^2,
$$
while the second one, together with Young's inequality, gives 
$U_4\leq U_{41}+U_{42}$, where
\begin{align*}
  U_{41} &= \frac{1}{4}\sum_\pm D(1\pm p)\sum_{\sigma\in\EE}\tau_{\sigma} 
	\min(n_{\pm,K}^{k,\delta},n_{\pm,K,\sigma}^{k,\delta})
	\big( \DD(\log n_{\pm}^{k,\delta}+V^k)_{K,\sigma}\big)^2, \\
  U_{42} &= \sum_\pm D(1\pm p)\sum_{\sigma\in\EE}\tau_{\sigma}
	\big(\max(n_{\pm,K}^{k,\delta},n_{\pm,K,\sigma}^{k,\delta})\big)^2 
	\frac{(\DD(\log n^D-\log (n^D+2\delta))_{K,\sigma})^2}{\min(n_{\pm,K}^{k,\delta},
	n_{\pm,K,\sigma}^{k,\delta})},\\
  &\leq 2D\frac{(M^0+\delta)^2}{\delta}\big|\log (n^D_\MM+\delta)
	-\log n^D_\MM\big|_{1,\MM}^2,
\end{align*}
since $\min(n_{\pm,K}^{k,\delta},n_{\pm,K,\sigma}^{k,\delta})\geq \delta$ 
for all $K\in\TT$ and $\sigma\in \EE_K$.
Applying Young's inequality again, we obtain $U_5\leq U_{51}+U_{52}$ with 
\begin{align*}
  U_{51} &= \frac{1}{4} \sum_{\pm} D(1\pm p) \sum_{\sigma\in\EE}\tau_{\sigma} 
	\min(n_{\pm,K}^{k,\delta},n_{\pm,K,\sigma}^{k,\delta})
	\big( \DD(\log n_{\pm}^{k,\delta}+V^k)_{K,\sigma}\big)^2, \\
  U_{52} &= \delta^2\sum_{\pm} D(1\pm p) \sum_{\sigma\in\EE}\tau_{\sigma}
	\frac{(\DD V_{K,\sigma}^k)^2}{\min(n_{\pm,K}^{k,\delta},
	n_{\pm,K,\sigma}^{k,\delta})}
  \leq 2D\delta \sum_{\sigma\in\EE}\tau_{\sigma}(\DD V_{K,\sigma}^k)^2
\end{align*}
and 
$$
  U_6 \leq 2D\delta\left(\big|\log (n^D_\MM+2\delta)-\log n^D_\MM\big|_{1,\MM}^2
	+ \sum_{\sigma\in\EE}\tau_{\sigma}(\DD V_{K,\sigma}^k)^2\right).
$$

Summarizing the above inequalities, we deduce that 
\begin{align}
  \frac{1}{\triangle t} & (E_\delta^k-E_\delta^{k-1})+\frac{1}{2} 
  \sum_{\pm} D(1\pm p) \sum_{\sigma\in\EE}\tau_{\sigma} 
	\min(n_{\pm,K}^{k,\delta},n_{\pm,K,\sigma}^{k,\delta})
	\big( \DD(\log n_{\pm}^{k,\delta}+V^k)_{K,\sigma}\big)^2 \label{thm2.li} \\
	&\leq 4D\delta\sum_{\sigma\in\EE}\tau_{\sigma}(\DD V_{K,\sigma}^k)^2
  + 2D\left(\delta+\frac{(M^0+\delta)^2}{\delta}\right)
	\big|\log (n^D_\MM+2\delta)-\log n^D_\MM\big|_{1,\MM}^2. \nonumber
\end{align}
On the one hand, the term $\sum_{\sigma\in\EE}\tau_{\sigma}(\DD V_{K,\sigma}^k)^2$ 
does not depend on $\delta$ and is bounded (this can be seen by using scheme 
\eqref{sch.v} and the $L^{\infty}$ bound on $n_{\pm,\TT}^k$).  
On the other hand, we rewrite 
\begin{align*}
  \big|\log (n^D_\MM+2\delta)-\log n^D_\MM\big|_{1,\MM}^2
	&= \left\vert \log \left(1+\frac{2\delta}{n_\MM^D}\right)\right\vert_{1,\MM}^2 \\
  &= \sum_{\sigma\in\EE}\tau_{\sigma} \left(\log 
	\left(1+\frac{2\delta}{n^D_{K,\sigma}}\right)
	- \log \left(1+\frac{2\delta}{n^D_{K}}\right)\right)^2.
\end{align*}
Employing the inequality $|\log y -\log x|\leq |y-x|/\min(x,y)$ for $x$, $y>0$, 
and the fact that $n^D\geq n_\ast>0$, we obtain 
$$
  \big|\log (n^D_\MM+2\delta)-\log n^D_\MM\big|_{1,\MM}^2
	\leq \frac{4\delta^2}{n_\ast^2}|n_\MM^D|_{1,\MM}^2.
$$
Thanks to hypothesis \eqref{hypo2}, $n^D\in H^1(\Omega)$, and Lemma 9.4 in 
\cite{EGH00}, we conclude that $\vert n_\MM^D\vert_{1,\MM}\leq 
K\Vert n^D\Vert_{H^1(\Omega)}$
with $K$ depending only on the regularity of the mesh $\MM$.
Therefore, the right-hand side in \eqref{thm2.li} tends to zero when $\delta\to 0$. 
Passing to the limit $\delta\to 0$ in \eqref{thm2.li} then leads to 
\eqref{ineq.thm.diss}. This concludes the proof of Theorem \ref{thm.diss}.

%%%%%%%%%%%%%%%%%%%%%%%%%%%%%%%%%%%%%%%%%%%%%%%%%%%%%%%%%%%%%%%%%%%%%%%%%%%%%%%

\section{Numerical simulations}\label{sec.numer}

\begin{figure}[ht]
  \centering
  \includegraphics[width=0.6\textwidth]{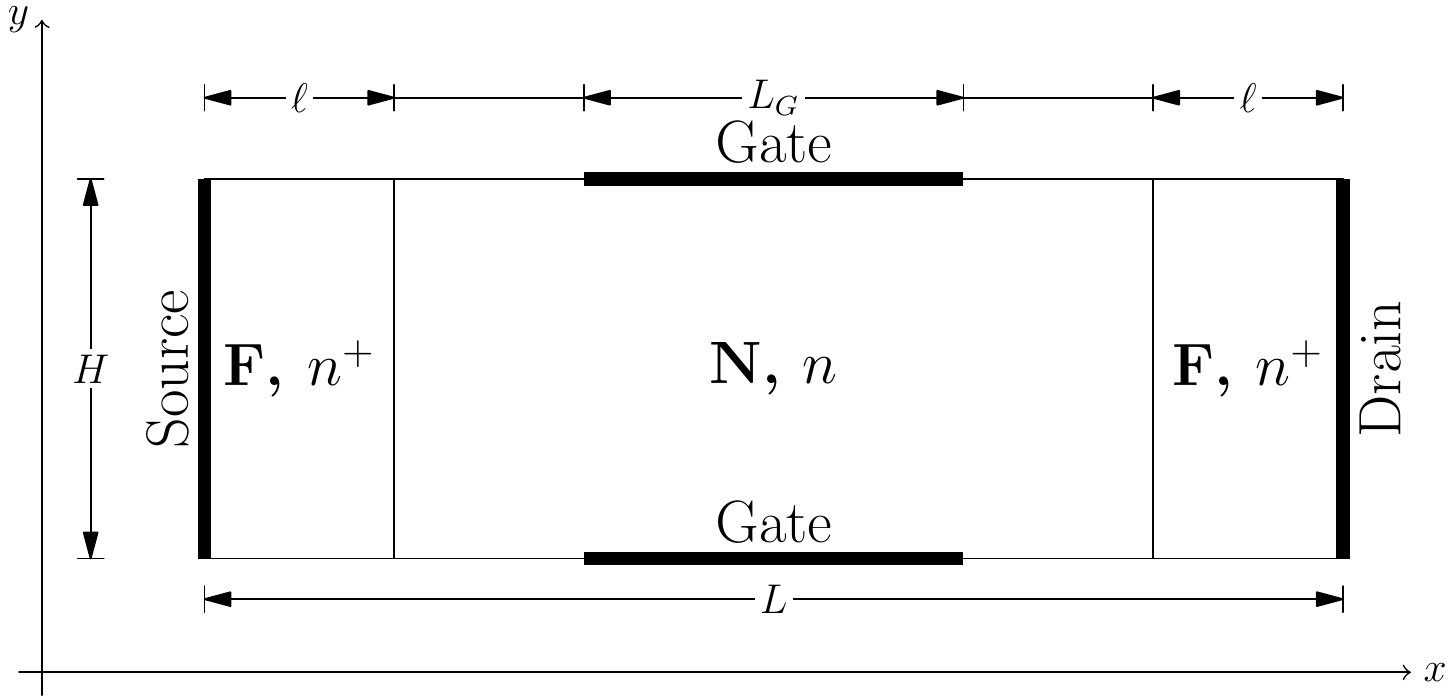} 
	 \caption{Geometry of a MESFET with ferromagnetic (F) source and drain regions
	and nonmagnetic (N) channel region.}
	\label{fig.MESFET}
\end{figure}

As an illustration of the numerical scheme, analyzed in the previous
sections, we present two-dimensional simulations of a simple double-gate
ferromagnetic MESFET (metal semiconductor field-effect transistor).
This device is composed of a semiconductor region which is sandwiched between
two ferromagnetic contact regions (see Figure \ref{fig.MESFET}).
The idea of such devices is that the source region plays the role of a spin 
polarizer. The non-zero spin-orbit interaction causes the electrons to precess during
the propagation through the middle channel region. At the drain contact,
only those electrons with spin aligned to the drain magnetization
can leave the channel and contribute to the current flow. Here, we focus
on the feasibility of our numerical scheme and the verification of the
properties of the numerical solution and less on the physical properties.
Therefore, the physical setting considered here is strongly simplified. 
In particular, we just modify the standard MESFET setup by allowing for 
ferromagnetic regions. For a more detailed modeling, we refer e.g.\ to \cite{LLN08}.

In the following, we describe the geometry of the device in the $(x,y)$ plane
(see Figure \ref{fig.MESFET}). The total length
is $L=0.6\,\mu$m and the height equals $H=0.2\,\mu$m. The source and drain regions
are highly doped with doping $C_+=3\cdot 10^{23}$m$^{-3}$. The doping in
the channel region is $C_0=10^{23}$m$^{-3}$. The length of the source and
drain regions are $\ell=0.1\,\mu$m. The gate contacts are attached at the
middle of the device with a length of $L_G=0.2\,\mu$m. 
The electrical parameters are: diffusion coefficient $D=10^{-3}\,$m$^2$s$^{-1}$,
relaxation time $\tau=1\,$ps, temperature $T=300\,$K, and relative
permittivity of silicon $\eps_r=11.7$. These parameters are similar to those
used in \cite{HJP04} (there is a small difference in the relaxation time value).
With these data, the scaled Debye length is 
$\lambda_D=(\eps_0\eps_r k_B T)/(qC_+ L^2)\approx 1.6\cdot 10^{-4}$, where
$\eps_0\approx 8.9\cdot 10^{-12}\,$Fm$^{-1}$ is the permittivity of the
vaccum, $q\approx 1.6\cdot 10^{-10}$\,C is the electron charge, and 
$k_B\approx 1.4\cdot 10^{-23}\,$m$^2$kgs$^{-2}$K$^{-1}$ is the Boltzmann
constant. 

The gate contact is considered as a Schottky contact with barrier potential
$V_S=0.8\,$V. The total voltage between source and gate is $V_G+V_S$, where
$V_G$ is the voltage applied at the gate. The density boundary value at the
gate contact is calclulated according to \cite[Formula (5.1-19)]{Sel84}, and
the potential of the closed state is taken from \cite{HJP04}. This gives
\begin{itemize}
  \item at the source: $n_0 = C_{+}$, $\vec n = 0$, potential: $0$\,V,
  \item at the drain: $n_0 = C_{+}$, $\vec n = 0$, potential: $V_D$,
  \item at the gate: \\
    open state: $n_0 = 3.9 \cdot 10^{11}\,$m$^{-3}$, $\vec n = 0$, potential: $V_S$, \\
    closed state: $n_0 = 3.2 \cdot 10^{9}\,$m$^{-3}$, $\vec n = 0$, 
		potential: $V_S + 1.2$\,V,
  \item for the other segments: homogeneous Neumann boundary conditions.
\end{itemize}

The magnetic field is caused by the local orientation of the electron spin
in the crystal and is predetermined by the ferromagnetic properties of the
material. We consider a constant magnetic field, oriented along the
$z$-axis (perpendicular to the device). The electron spin may be also
changed under the influence of the spin current, but we do not consider
this effect here. In our model, $\vec{m}$ corresponds to the direction of
the local magnetic field, and the parameter $\gamma$ describes the intensity
of the spin precession around this field. We choose $\vec{m}=0$ in the
channel region and 
$$ 
  \vec m = \left\{
  \begin{aligned}
    &(0, 0, 1) \quad\mbox{for } x < L/3 \text{ or } x \ge 2L/3, \\
    &(0, 0, 0) \quad\mbox{for } L/3 \le x < 2L/3.
  \end{aligned}\right. 
$$
The value for $\gamma$ is taken from \cite{PoNe11}, i.e.\ $\gamma=\hbar/\tau$,
with $\hbar$ being the reduced Planck constant. The spin polarization is
nonzero only in the highly doped source and drain regions, and we take $p=0.9$. 

For the numerical discretization, we have chosen an admissible triangular mesh.
Equations \eqref{1.n0}-\eqref{1.V2} are approximated by scheme 
\eqref{sch.n0}-\eqref{sch.J}, with the corresponding boundary conditions.
The nonlinear system is solved at each time step by Newton's method.
The time step size is $\triangle t=0.05$. The computations are continued
until a steady state is reached or, more precisely,
until the difference of the solutions at two consecutive time steps in the 
$\ell^2$ norm falls below a threshold (typically, $10^{-5}$).

\begin{figure}[hT]
  \centering
  \includegraphics[width=0.49\textwidth]{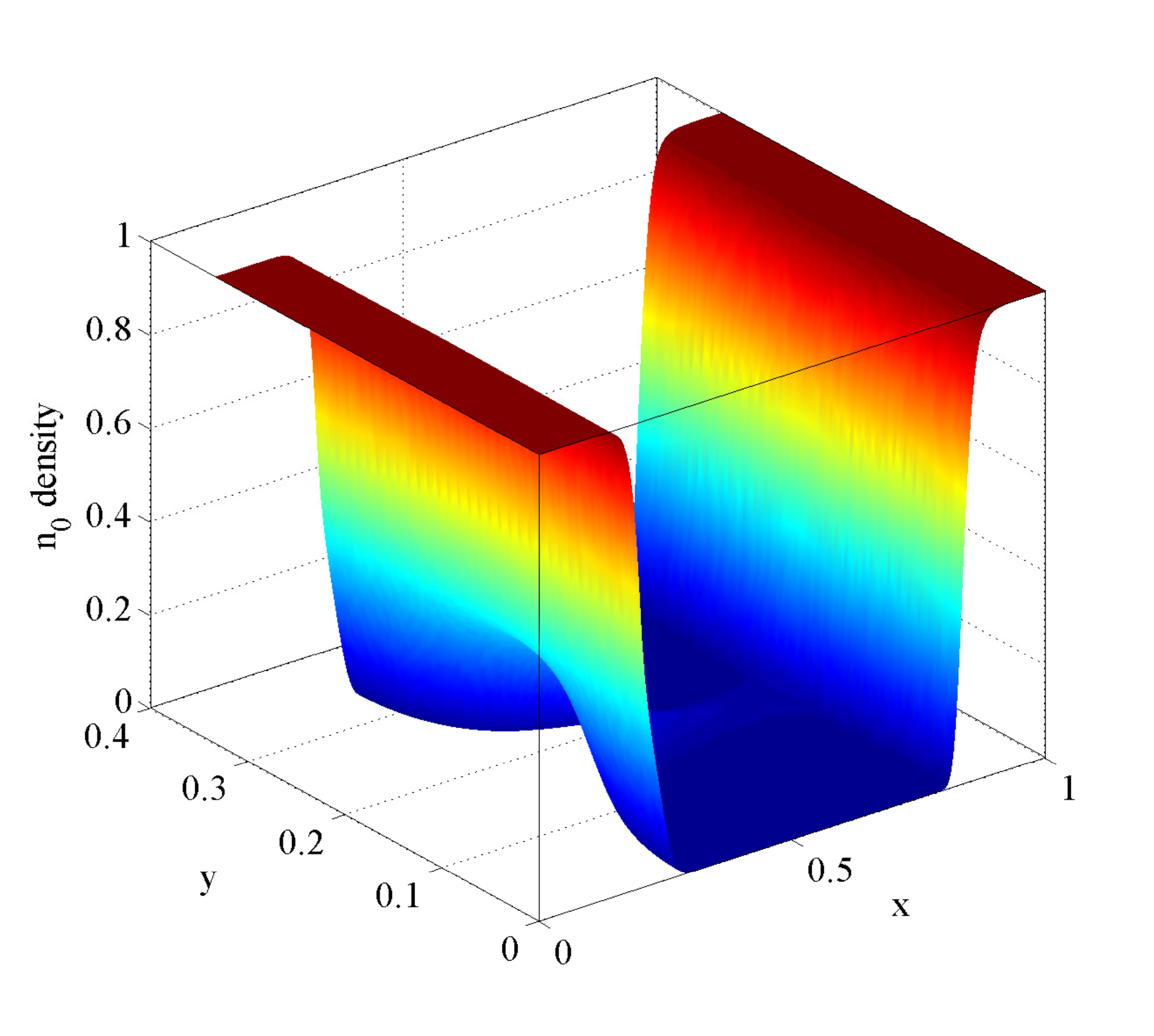}
  \includegraphics[width=0.49\textwidth]{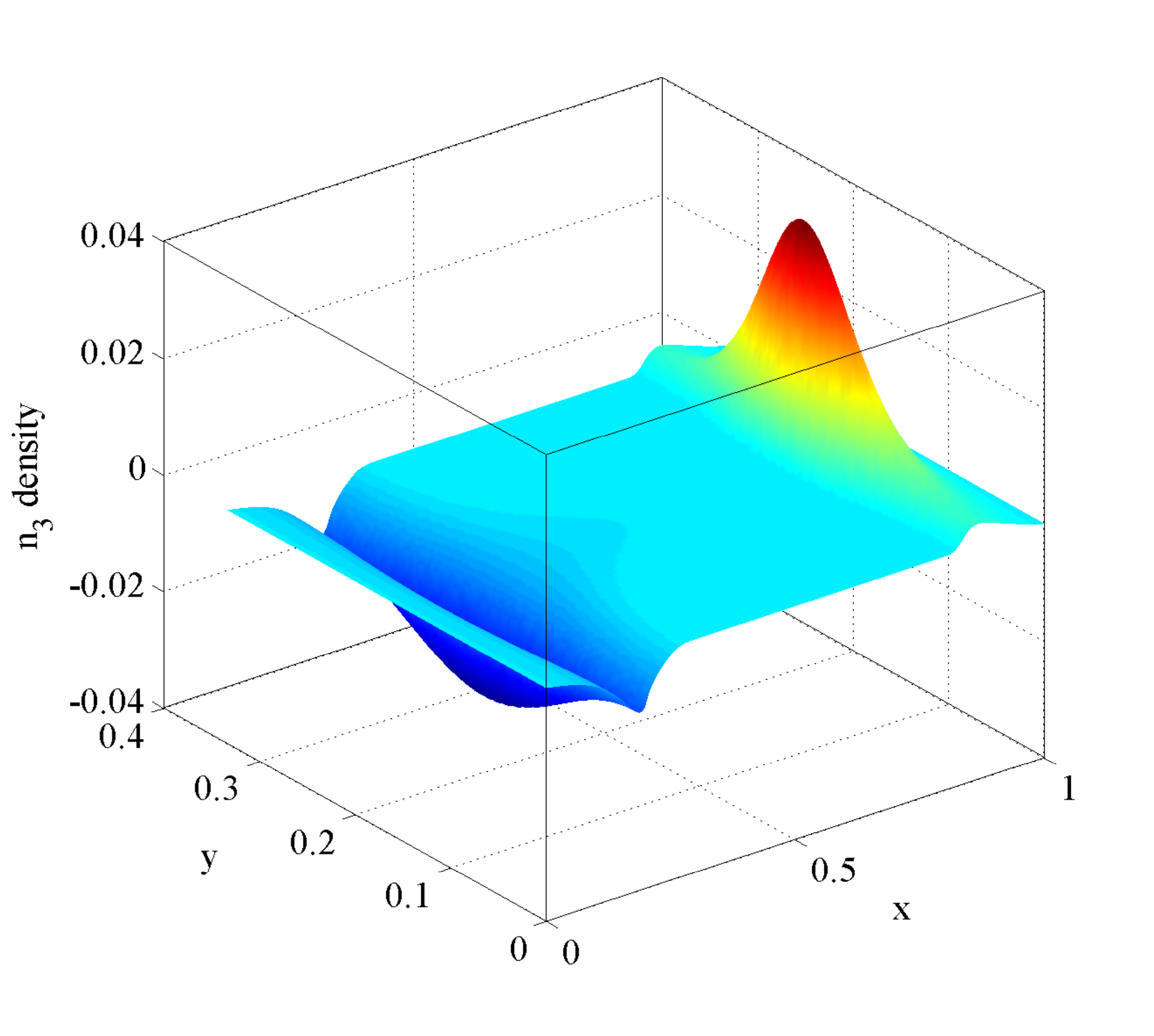} 
	\caption{Scaled stationary charge density (left) and spin density 
	$n_3$ (right) in an open-state MESFET with $V_D=-2$\,V and $V_G=0$\,V.}
	\label{fig.steady.o}
\end{figure}

\begin{figure}[ht]
  \centering
  \includegraphics[width=0.49\textwidth]{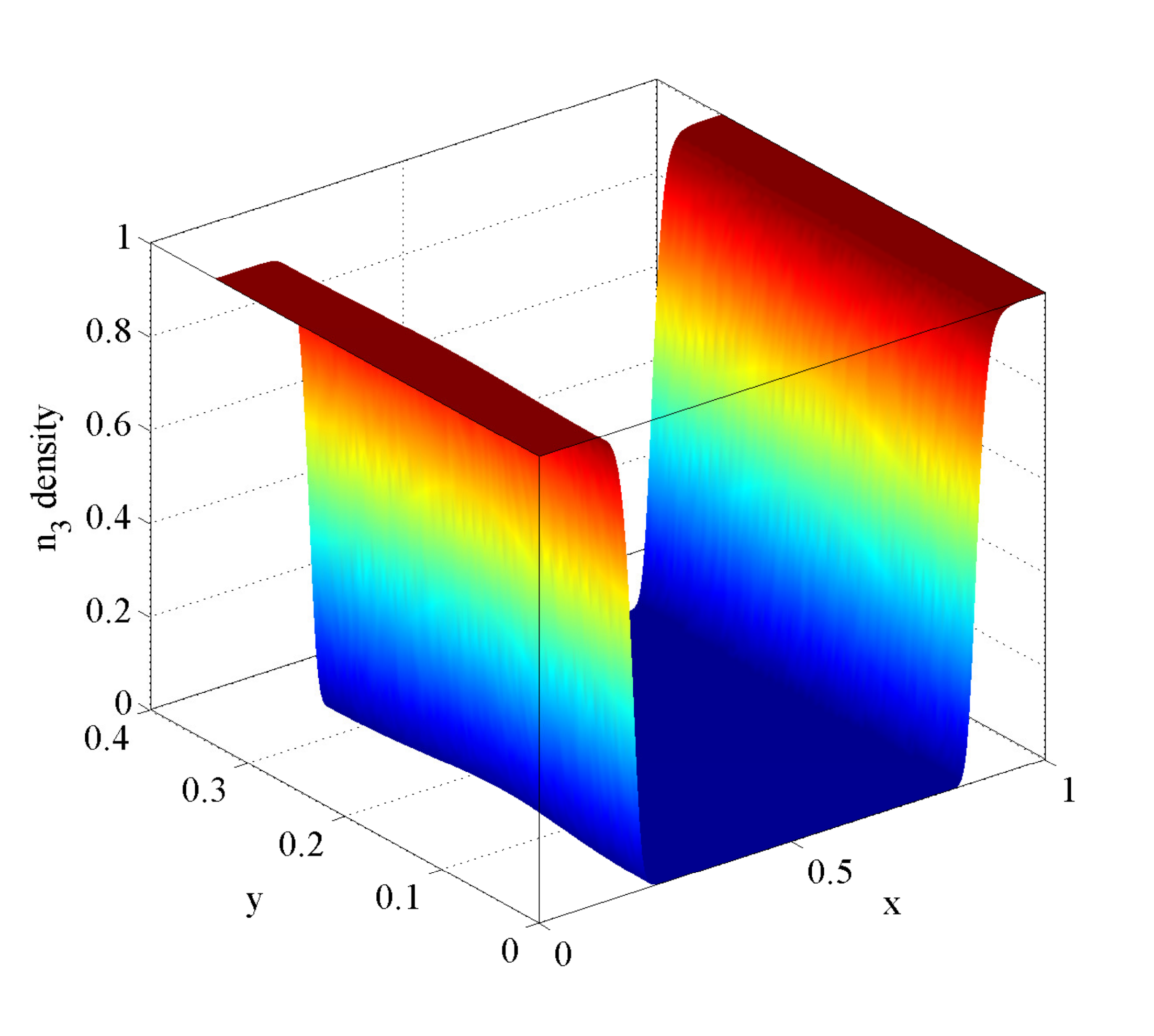}
	\caption{Scaled stationary charge density in a closed-state MESFET 
	with $V_D=-2$\,V and $V_G=1.2$\,V.}
	\label{fig.steady.c}
\end{figure}

Figure \ref{fig.steady.o} illustrates the scaled steady-state charge density $n_0$
and the spin density $n_3$ (note that $n_1=n_2=0$) in the open state.
The densities are scaled by the doping concentration $C_+$, the spatial
variable by the device length $L$.
Compared to the closed state in Figure \ref{fig.steady.c}, the charge
density is rather large in the channel region, which can be also observed
in standard MESFET devices. The charge current density in the closed state 
is by the factor $10^5\ldots 10^6$ smaller than in the open state. The spin
density is (almost) zero in the closed state. 

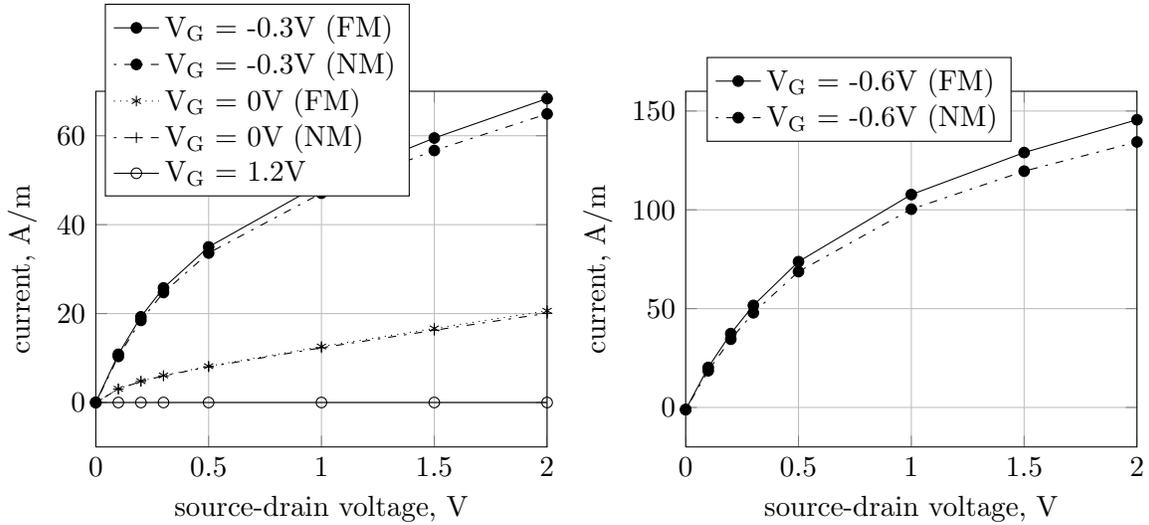
\begin{figure}[ht]
  \centering\small
  \begin{tikzpicture}

\begin{axis}[%
width=\vacwidthall,
height=0.788709677419355\vacwidthall,
scale only axis,
xmin=0,
xmax=2,
xlabel={source-drain voltage, V},
xmajorgrids,
ymin=-10,
ymax=70,
ylabel={current, A/m},
ymajorgrids,
legend style={at={(0.020002380952378,0.704365079365079)},anchor=south west,draw=black,fill=white,legend cell align=left}
]
\addplot [
color=black,
solid,
mark=*,
mark options={solid}
]
table[row sep=crcr]{
0 -0.000430458\\
0.1 10.84657766\\
0.2 19.27092149\\
0.3 25.76972488\\
0.5 34.99760588\\
1 49.21242577\\
1.5 59.48435045\\
2 68.35882927\\
};
\addlegendentry{$\text{V}_\text{G}\text{ = -0.3V (FM)}$};

\addplot [
color=black,
dash pattern=on 1pt off 3pt on 3pt off 3pt,
mark=*,
mark options={solid}
]
table[row sep=crcr]{
0 1.40408e-05\\
0.1 10.35271392\\
0.2 18.46073563\\
0.3 24.73773616\\
0.5 33.61666628\\
1 47.08264736\\
1.5 56.67877171\\
2 64.90436064\\
};
\addlegendentry{$\text{V}_\text{G}\text{ = -0.3V (NM)}$};

\addplot [
color=black,
dotted,
mark=asterisk,
mark options={solid}
]
table[row sep=crcr]{
0 5.02756e-05\\
0.1 3.064894873\\
0.2 4.803559104\\
0.3 6.073921527\\
0.5 8.165950618\\
1 12.58641948\\
1.5 16.66358661\\
2 20.62955454\\
};
\addlegendentry{$\text{V}_\text{G}\text{ = 0V (FM)}$};

\addplot [
color=black,
dash pattern=on 1pt off 3pt on 3pt off 3pt,
mark=+,
mark options={solid}
]
table[row sep=crcr]{
0 -0.000109669\\
0.1 3.020097689\\
0.2 4.733509531\\
0.3 5.979416849\\
0.5 8.022361089\\
1 12.31424711\\
1.5 16.24912142\\
2 20.05864361\\
};
\addlegendentry{$\text{V}_\text{G}\text{ = 0V (NM)}$};

\addplot [
color=black,
solid,
mark=o,
mark options={solid}
]
table[row sep=crcr]{
0 8.38346e-05\\
0.1 1.06044e-05\\
0.2 3.67859e-05\\
0.3 1.35857e-05\\
0.5 8.9063e-06\\
1 1.50765e-05\\
1.5 1.9031e-05\\
2 1.11187e-05\\
};
\addlegendentry{$\text{V}_\text{G}\text{ = 1.2V}$};

\end{axis}
\end{tikzpicture}
	\begin{tikzpicture}

\begin{axis}[%
width=\vacwidthneg,
height=0.788709677419355\vacwidthneg,
scale only axis,
xmin=0,
xmax=2,
xlabel={source-drain voltage, V},
xmajorgrids,
ymin=-20,
ymax=160,
ylabel={current, A/m},
ymajorgrids,
legend style={at={(0.046547619047617,0.857936507936508)},anchor=south west,draw=black,fill=white,legend cell align=left}
]
\addplot [
color=black,
solid,
mark=*,
mark options={solid}
]
table[row sep=crcr]{
0 -1.084915615\\
0.1 20.23339396\\
0.2 37.37160164\\
0.3 51.64175714\\
0.5 73.77184943\\
1 107.7213514\\
1.5 128.9915867\\
2 145.5893486\\
};
\addlegendentry{$\text{V}_\text{G}\text{ = -0.6V (FM)}$};

\addplot [
color=black,
dash pattern=on 1pt off 3pt on 3pt off 3pt,
mark=*,
mark options={solid}
]
table[row sep=crcr]{
0 -1.076585654\\
0.1 18.56441097\\
0.2 34.5129809\\
0.3 47.89171138\\
0.5 68.73167793\\
1 100.3221636\\
1.5 119.5747032\\
2 134.3488976\\
};
\addlegendentry{$\text{V}_\text{G}\text{ = -0.6V (NM)}$};

\end{axis}
\end{tikzpicture}
	\caption{Current-voltage characteristics for the ferromagnetic (FM) and
	standard (NM) MESFET for various gate voltages $V_G$. For convenience, 
	the source-drain voltages are given by their absolute values.}
	\label{fig.VAC}
\end{figure}

Current-voltage characteristics for MESFETs with and without ferromagnetic regions
are shown in Figure \ref{fig.VAC}. We observe that in the open state
(nonpositive gate potentials), the current densities
in the ferromagnetic MESFET are slightly larger than in the standard device,
which allows for an improved device performance. When the transistor is closed
($V_G=1.2$\,V), the current densities are (almost) zero for both transistor types.

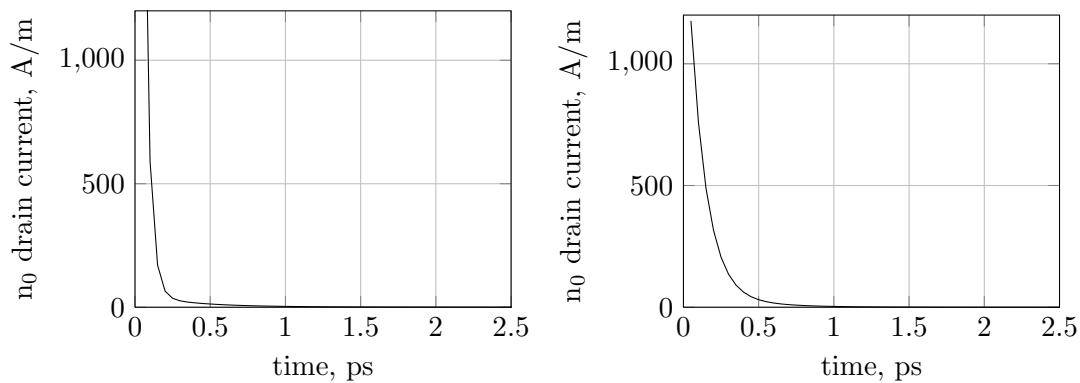
\begin{figure}[ht]
  \centering\small
  \begin{tikzpicture}

\begin{axis}[%
width=\optoclwidth,
height=0.788709677419355\optoclwidth,
scale only axis,
xmin=0,
xmax=2.5,
xlabel={time, ps},
xmajorgrids,
ymin=0,
ymax=1200,
ylabel={$\text{n}_\text{0}\text{ drain current, A/m}$},
ymajorgrids
]
\addplot [
color=black,
solid,
forget plot
]
table[row sep=crcr]{
0.05 2332.09080426599\\
0.1 591.856175925042\\
0.15 170.756369165644\\
0.2 65.4004641228111\\
0.25 36.147789801966\\
0.3 25.9291224163685\\
0.35 20.9306523173709\\
0.4 17.6615756196427\\
0.45 15.1562599859519\\
0.5 13.1001454150074\\
0.55 11.3638269828906\\
0.6 9.87929875514613\\
0.65 8.60314072297999\\
0.7 7.5037266658906\\
0.75 6.55601335122752\\
0.8 5.73914316077309\\
0.85 5.03524721001206\\
0.9 4.42880484773147\\
0.95 3.90626842014225\\
1 3.45581371034055\\
1.05 3.06714889358425\\
1.1 2.73135090794948\\
1.15 2.44071644923144\\
1.2 2.1886228437998\\
1.25 1.96939828755683\\
1.3 1.77820166347193\\
1.35 1.61091282690333\\
1.4 1.46403359601923\\
1.45 1.33459947776083\\
1.5 1.22010182907413\\
1.55 1.11841981913856\\
1.6 1.02776145445573\\
1.65 0.946612731280143\\
1.7 0.87369409392135\\
1.75 0.807923399734375\\
1.8 0.748384421734709\\
1.85 0.694300191602473\\
1.9 0.645010618653008\\
1.95 0.599953602044675\\
2 0.558649125440468\\
2.05 0.520686179364096\\
2.1 0.485711439065709\\
2.15 0.453420220755477\\
2.2 0.42354857359359\\
2.25 0.395867036093646\\
2.3 0.370175078725017\\
2.35 0.346296736689427\\
2.4 0.324076837972543\\
2.45 0.303378016981831\\
2.5 0.284077922467837\\
2.55 0.266067189755518\\
2.6 0.24924760556918\\
2.65 0.233530642024696\\
2.7 0.218836047780713\\
2.75 0.205091011882631\\
2.8 0.192228958296263\\
2.85 0.180189080660938\\
2.9 0.168915467618303\\
2.95 0.158356756240474\\
3 0.148465515649131\\
3.05 0.13919784781192\\
3.1 0.130513144943095\\
3.15 0.122373715193499\\
3.2 0.114744504646517\\
3.25 0.107592970359632\\
3.3 0.100888450074623\\
3.35 0.0946027638269356\\
3.4 0.0887096003738418\\
3.45 0.0831840640728774\\
3.5 0.0780030728458739\\
3.55 0.073145144379607\\
3.6 0.0685897694843713\\
3.65 0.0643181433685257\\
3.7 0.0603126859515409\\
3.75 0.056556633795995\\
3.8 0.0530345541062087\\
3.85 0.0497316995494018\\
3.9 0.0466345931611359\\
3.95 0.0437301874597598\\
4 0.041006651699373\\
4.05 0.0384524159606524\\
4.1 0.0360577404340469\\
4.15 0.0338118433126544\\
4.2 0.0317057592850297\\
4.25 0.0297311520660812\\
4.3 0.0278790293504892\\
4.35 0.0261426668279824\\
4.4 0.024513869592473\\
4.45 0.0229869000613332\\
4.5 0.0215554543357607\\
4.55 0.0202121115988159\\
4.6 0.0189533417545972\\
4.65 0.0177723776015181\\
4.7 0.0166652559404738\\
4.75 0.0156268253005056\\
4.8 0.0146537899161649\\
4.85 0.0137402522958551\\
4.9 0.0128846280981769\\
4.95 0.0120815918029891\\
5 0.0113288386509857\\
5.05 0.0106229758818314\\
5.1 0.00996127678434187\\
5.15 0.00934109799279176\\
5.2 0.00875918887257079\\
5.25 0.00821276110097288\\
5.3 0.00770166751260002\\
5.35 0.0072216135827752\\
5.4 0.00677126405662547\\
5.45 0.006349433003815\\
5.5 0.00595366998515482\\
5.55 0.00558344797705454\\
5.6 0.00523503005258545\\
5.65 0.00490880877319985\\
5.7 0.00460303978250252\\
5.75 0.00431646329058862\\
5.8 0.00404688879060847\\
5.85 0.00379524389264856\\
5.9 0.00355903736930526\\
5.95 0.00333721775892326\\
6 0.00312888031791693\\
6.05 0.00293370152102009\\
6.1 0.00275131620102862\\
6.15 0.00257961431585528\\
6.2 0.00241874936445119\\
6.25 0.0022679156030768\\
6.3 0.00212643001850464\\
6.35 0.00199391783708288\\
6.4 0.00186994142014746\\
6.45 0.00175348925905925\\
6.5 0.00164386442145675\\
6.55 0.00154193978452307\\
6.6 0.00144558416208147\\
6.65 0.00135576450597429\\
6.7 0.00127092762790661\\
6.75 0.00119144665016022\\
6.8 0.00111751352249617\\
6.85 0.00104781228570866\\
6.9 0.00098257941910394\\
6.95 0.000920866041496687\\
7 0.000864194699107805\\
7.05 0.000810172662688581\\
7.1 0.000759756843721323\\
7.15 0.000711658781459873\\
7.2 0.000668272339729718\\
7.25 0.000625702083615377\\
7.3 0.000587234868961779\\
7.35 0.000550232517191983\\
7.4 0.000516034301393862\\
7.45 0.000484006958770116\\
7.5 0.000453590531933984\\
7.55 0.000425486919428873\\
7.6 0.000399229459506367\\
7.65 0.000374208081597644\\
7.7 0.000350855806347364\\
7.75 0.000329002554600291\\
7.8 0.000308226329232866\\
7.85 0.000289032764974596\\
7.9 0.00027061501193998\\
7.95 0.000254535264237767\\
8 0.000238552200127373\\
8.05 0.000223890590039429\\
8.1 0.00021010327489398\\
8.15 0.000197045119991701\\
8.2 0.000184186236601736\\
8.25 0.000172124002095874\\
8.3 0.000161442441605513\\
8.35 0.000151503090499948\\
8.4 0.000142638472287243\\
8.45 0.000133188460341693\\
8.5 0.000125306882655561\\
8.55 0.000117505245988201\\
8.6 0.000110033324314239\\
8.65 0.000103682803257471\\
8.7 9.72979903438455e-05\\
8.75 9.06997550343471e-05\\
8.8 8.65314628888655e-05\\
8.85 7.99667152738773e-05\\
8.9 7.46158125500332e-05\\
8.95 7.03860839018517e-05\\
9 6.57649237878249e-05\\
9.05 6.21974089611471e-05\\
9.1 5.71650625665116e-05\\
9.15 5.34553062082255e-05\\
9.2 5.08682111798711e-05\\
9.25 4.81112653740334e-05\\
9.3 4.51167286966724e-05\\
9.35 4.10711107197984e-05\\
9.4 3.94967924904346e-05\\
9.45 3.72864562412331e-05\\
9.5 3.49417025741484e-05\\
9.55 3.20774467757213e-05\\
9.6 3.04129210771109e-05\\
9.65 2.91565121570687e-05\\
9.7 2.66596078424495e-05\\
9.75 2.51511119645615e-05\\
9.8 2.34535567295542e-05\\
9.85 2.20043610024964e-05\\
9.9 2.05589641629961e-05\\
9.95 1.92022840966215e-05\\
10 1.82858688790419e-05\\
10.05 1.71522662925952e-05\\
10.1 1.60803110435278e-05\\
10.15 1.48340236246157e-05\\
10.2 1.37131889953605e-05\\
10.25 1.31871514648002e-05\\
10.3 1.22139408000127e-05\\
10.35 1.14079345055642e-05\\
10.4 1.09541719765158e-05\\
10.45 1.02284988966463e-05\\
10.5 9.66981281785216e-06\\
10.55 9.07606210925782e-06\\
10.6 7.43238083005666e-06\\
10.65 7.81870571143571e-06\\
10.7 7.58687800371511e-06\\
10.75 6.22514062566601e-06\\
};
\end{axis}
\end{tikzpicture} 
	\begin{tikzpicture}

\begin{axis}[%
width=\optoclwidth,
height=0.777442396313364\optoclwidth,
scale only axis,
xmin=0,
xmax=2.5,
xlabel={time, ps},
xmajorgrids,
ymin=0,
ymax=1200,
ylabel={$\text{n}_\text{0}\text{ drain current, A/m}$},
ymajorgrids
]
\addplot [
color=black,
solid,
forget plot
]
table[row sep=crcr]{
0.05 1175.96902914721\\
0.1 758.105665220339\\
0.15 488.570980092223\\
0.2 315.920128391926\\
0.25 205.745108776884\\
0.3 135.509907830659\\
0.35 90.6586307779441\\
0.4 61.8801328890167\\
0.45 43.259642025569\\
0.5 31.0593611287732\\
0.55 22.9264068796962\\
0.6 17.3837674422247\\
0.65 13.505683267122\\
0.7 10.711896055911\\
0.75 8.637912669989\\
0.8 7.0535835204267\\
0.85 5.8122309348175\\
0.9 4.81900334072591\\
0.95 4.01127256812377\\
1 3.34653914050502\\
1.05 2.79499477195631\\
1.1 2.33495853119573\\
1.15 1.95007633942722\\
1.2 1.62759642723688\\
1.25 1.35729785742911\\
1.3 1.13081372604074\\
1.35 0.941192221099713\\
1.4 0.78260136569603\\
1.45 0.650121059536427\\
1.5 0.539589178468558\\
1.55 0.447482261768643\\
1.6 0.370818981988076\\
1.65 0.307080299242052\\
1.7 0.254141186973302\\
1.75 0.210212790069905\\
1.8 0.173792216466323\\
1.85 0.143619209341664\\
1.9 0.11863896385726\\
1.95 0.0979704080250526\\
2 0.0808783997155295\\
2.05 0.0667507589595232\\
2.1 0.0550783058388343\\
2.15 0.0454377890937448\\
2.2 0.0374781210538006\\
2.25 0.0309080930954683\\
2.3 0.025486352026304\\
2.35 0.0210130512833841\\
2.4 0.0173236034549564\\
2.45 0.0142804377605744\\
2.5 0.0117709337202303\\
2.55 0.00970188655698814\\
2.6 0.00799586903979504\\
2.65 0.00658984754499745\\
2.7 0.0054304268283694\\
2.75 0.00447500280125421\\
2.8 0.00368764156353003\\
2.85 0.00303876693832929\\
2.9 0.00250407227434718\\
2.95 0.00206314527828811\\
3 0.00169996183537046\\
3.05 0.0014007604553536\\
3.1 0.00115397096176735\\
3.15 0.000950838207497824\\
3.2 0.000783550911642915\\
3.25 0.00064550033318161\\
3.3 0.000531825087512375\\
3.35 0.000438153665134183\\
3.4 0.000360990546866265\\
3.45 0.000297283331150472\\
3.5 0.000245009906095548\\
3.55 0.000201908744362778\\
3.6 0.000166256336217232\\
3.65 0.000137062328339577\\
3.7 0.000112979368606604\\
3.75 9.31121818842847e-05\\
3.8 7.68341897834665e-05\\
3.85 6.31454129781325e-05\\
3.9 5.18538010393976e-05\\
3.95 4.27561897617923e-05\\
4 3.52627939113405e-05\\
4.05 2.90903685405559e-05\\
4.1 2.39760969476938e-05\\
4.15 1.96552256602872e-05\\
4.2 1.62931799784474e-05\\
4.25 1.34639421477822e-05\\
4.3 1.10304627831133e-05\\
};
\end{axis}
\end{tikzpicture}
	\caption{Change of the electron drain current in the ferromagnetic (left)
	and standard (right) MESFET, switching from open to closed state.}
	\label{fig.trans1}
\end{figure}

In Figure \ref{fig.trans1} (left), we present the transient behavior of the
charge density when switching from the open to the closed state ($V_D=-2$\,V). 
The current values stabilize after about 1\,ps.
This justifies to define the numerical solution after 12\,ps as the ``steady-state
solution''. We compare these values with those computed from a standard MESFET;
see Figure \ref{fig.trans1} (right). The stabilization in the ferromagnetic case is
slightly faster which allows for faster devices.

\begin{figure}[ht]
  \centering\small
\begin{tikzpicture}

\begin{axis}[%
width=\entropywidth,
height=0.688709677419355\entropywidth,
scale only axis,
xmin=0,
xmax=25,
xlabel={time, ps},
xmajorgrids,
ymode=log,
ymin=1e-20,
ymax=1,
yminorticks=true,
ylabel={free energy},
ymajorgrids,
yminorgrids
]
\addplot [
color=black,
solid,
forget plot
]
table[row sep=crcr]{
0.05 0.544387214979517\\
0.1 0.455023133569535\\
0.15 0.387719213644896\\
0.2 0.332512171403657\\
0.25 0.286309805497352\\
0.3 0.247307872786722\\
0.35 0.214211068058735\\
0.4 0.186017443300806\\
0.45 0.16192461920454\\
0.5 0.141278126364036\\
0.55 0.123538526656914\\
0.6 0.108258322422383\\
0.65 0.0950646484169352\\
0.7 0.0836457698870254\\
0.75 0.0737403053893174\\
0.8 0.0651285131435685\\
0.85 0.0576251929761615\\
0.9 0.051073876513315\\
0.95 0.0453420550239999\\
1 0.0403172481170013\\
1.05 0.0359037566505798\\
1.1 0.0320199742216628\\
1.15 0.0285961559697875\\
1.2 0.0255725627604665\\
1.25 0.0228979142350015\\
1.3 0.0205280965743124\\
1.35 0.0184250807696877\\
1.4 0.016556015224342\\
1.45 0.0148924630156171\\
1.5 0.0134097594333755\\
1.55 0.0120864697152497\\
1.6 0.0109039304136541\\
1.65 0.00984586070373512\\
1.7 0.00889803229662152\\
1.75 0.0080479885554212\\
1.8 0.00728480500077934\\
1.85 0.00659888470159536\\
1.9 0.00598178312619573\\
1.95 0.00542605792126224\\
2 0.00492513982414886\\
2.05 0.00447322152633378\\
2.1 0.00406516181410899\\
2.15 0.00369640273556342\\
2.2 0.0033628978954052\\
2.25 0.0030610502735342\\
2.3 0.00278765820935057\\
2.35 0.00253986840017572\\
2.4 0.00231513493523957\\
2.45 0.00211118353233564\\
2.5 0.00192598026691842\\
2.55 0.00175770418693824\\
2.6 0.00160472329425004\\
2.65 0.00146557344756\\
2.7 0.00133893980477304\\
2.75 0.00122364047606042\\
2.8 0.00111861210447237\\
2.85 0.00102289712973681\\
2.9 0.000935632524033288\\
2.95 0.000856039816901402\\
3 0.000783416250753953\\
3.05 0.000717126929336173\\
3.1 0.000656597839416666\\
3.15 0.000601309641446963\\
3.2 0.00055079213825612\\
3.25 0.000504619342357873\\
3.3 0.00046240507240535\\
3.35 0.000423799017959611\\
3.4 0.000388483219219592\\
3.45 0.00035616891486853\\
3.5 0.000326593716847593\\
3.55 0.000299519075789611\\
3.6 0.000274728005181322\\
3.65 0.000252023036013526\\
3.7 0.000231224377034877\\
3.75 0.0002121682585574\\
3.8 0.000194705440303617\\
3.85 0.000178699865994294\\
3.9 0.000164027449326304\\
3.95 0.000150574977702222\\
4 0.000138239121583344\\
4.05 0.000126925538669857\\
4.1 0.00011654806328736\\
4.15 0.000107027972399002\\
4.2 9.82933205805921e-05\\
4.25 9.02783371142997e-05\\
4.3 8.29228790745407e-05\\
4.35 7.61719349258672e-05\\
4.4 6.99751737174797e-05\\
4.45 6.42865354712089e-05\\
4.5 5.90638588058184e-05\\
4.55 5.42685422469924e-05\\
4.6 4.98652360276658e-05\\
4.65 4.58215615057249e-05\\
4.7 4.21078556110397e-05\\
4.75 3.86969379912847e-05\\
4.8 3.55638987456082e-05\\
4.85 3.26859048693431e-05\\
4.9 3.00420236673272e-05\\
4.95 2.76130616090577e-05\\
5 2.53814172165169e-05\\
5.05 2.33309467187314e-05\\
5.1 2.14468413252608e-05\\
5.15 1.97155150802233e-05\\
5.2 1.81245023530397e-05\\
5.25 1.66623641125474e-05\\
5.3 1.53186022074317e-05\\
5.35 1.40835809497591e-05\\
5.4 1.29484553605441e-05\\
5.45 1.19051054968307e-05\\
5.5 1.09460763305139e-05\\
5.55 1.00645226962108e-05\\
5.6 9.25415887548395e-06\\
5.65 8.50921240585814e-06\\
5.7 7.82438177567342e-06\\
5.75 7.19479763787933e-06\\
5.8 6.61598727475309e-06\\
5.85 6.08384201797155e-06\\
5.9 5.59458737906624e-06\\
5.95 5.14475566002751e-06\\
6 4.73116083452758e-06\\
6.05 4.35087550811211e-06\\
6.1 4.0012097828418e-06\\
6.15 3.67969186667952e-06\\
6.2 3.38405027913313e-06\\
6.25 3.11219752598691e-06\\
6.3 2.8622151120767e-06\\
6.35 2.63233978728928e-06\\
6.4 2.4209509199405e-06\\
6.45 2.22655890550062e-06\\
6.5 2.04779452400328e-06\\
6.55 1.88339916757174e-06\\
6.6 1.73221586716682e-06\\
6.65 1.59318105146175e-06\\
6.7 1.46531697775683e-06\\
6.75 1.3477247806663e-06\\
6.8 1.23957808545345e-06\\
6.85 1.14011714414818e-06\\
6.9 1.04864344381706e-06\\
6.95 9.64514757452837e-07\\
7 8.87140594748444e-07\\
7.05 8.15978022674056e-07\\
7.1 7.50527824661402e-07\\
7.15 6.9033097167872e-07\\
7.2 6.34965379130316e-07\\
7.25 5.84042926314835e-07\\
7.3 5.37206717471136e-07\\
7.35 4.94128564147642e-07\\
7.4 4.54506671387109e-07\\
7.45 4.1806351058132e-07\\
7.5 3.84543865250015e-07\\
7.55 3.53713033707878e-07\\
7.6 3.25355176794039e-07\\
7.65 2.99271799905877e-07\\
7.7 2.75280357171345e-07\\
7.75 2.53212967464309e-07\\
7.8 2.32915235040838e-07\\
7.85 2.14245164556634e-07\\
7.9 1.97072164142371e-07\\
7.95 1.8127612802742e-07\\
8 1.66746593666191e-07\\
8.05 1.53381966770773e-07\\
8.1 1.41088808311172e-07\\
8.15 1.29781179276697e-07\\
8.2 1.19380038075343e-07\\
8.25 1.09812686347468e-07\\
8.3 1.01012259592684e-07\\
8.35 9.29172581543785e-08\\
8.4 8.54711170909128e-08\\
8.45 7.86218094379521e-08\\
8.5 7.23214818769817e-08\\
8.55 6.65261197687516e-08\\
8.6 6.1195238975424e-08\\
8.65 5.62916026482808e-08\\
8.7 5.17809597870163e-08\\
8.75 4.76318070748375e-08\\
8.8 4.38151673840618e-08\\
8.85 4.03043869160406e-08\\
8.9 3.70749492988531e-08\\
8.95 3.41043039700229e-08\\
9 3.13717087415436e-08\\
9.05 2.88580837059003e-08\\
9.1 2.65458795481416e-08\\
9.15 2.44189529543545e-08\\
9.2 2.24624554015008e-08\\
9.25 2.06627285492669e-08\\
9.3 1.90072088966154e-08\\
9.35 1.74843395327193e-08\\
9.4 1.60834908722544e-08\\
9.45 1.47948845206445e-08\\
9.5 1.36095259916772e-08\\
9.55 1.25191421971247e-08\\
9.6 1.15161223508892e-08\\
9.65 1.05934659953735e-08\\
9.7 9.74473382115169e-09\\
9.75 8.96400248563934e-09\\
9.8 8.24582305064178e-09\\
9.85 7.58518378401974e-09\\
9.9 6.97747469808907e-09\\
9.95 6.41845448075129e-09\\
10 5.90422101075119e-09\\
10.05 5.43118705504468e-09\\
10.1 4.9960513419966e-09\\
10.15 4.59577743139934e-09\\
10.2 4.22757268100945e-09\\
10.25 3.88886718795469e-09\\
10.3 3.57729754057399e-09\\
10.35 3.29068848949513e-09\\
10.4 3.02704223552953e-09\\
10.45 2.78451822025641e-09\\
10.5 2.56144561331182e-09\\
10.55 2.35620101819219e-09\\
10.6 2.16742191600268e-09\\
10.65 1.99376694694677e-09\\
10.7 1.83402485979766e-09\\
10.75 1.68709576358885e-09\\
10.8 1.55190829439563e-09\\
10.85 1.42756565475519e-09\\
10.9 1.3131857024594e-09\\
10.95 1.20796823571727e-09\\
11 1.11118308565261e-09\\
11.05 1.0221478415007e-09\\
11.1 9.40247984817899e-10\\
11.15 8.64910429727192e-10\\
11.2 7.95607601468758e-10\\
11.25 7.3185659162588e-10\\
11.3 6.73213974874598e-10\\
11.35 6.19271801080999e-10\\
11.4 5.69647478559991e-10\\
11.45 5.24000782423988e-10\\
11.5 4.82034198061747e-10\\
11.55 4.43391988360929e-10\\
11.6 4.07856850098957e-10\\
11.65 3.75182713610127e-10\\
11.7 3.45106748442463e-10\\
11.75 3.17450463552444e-10\\
11.8 2.92017350823388e-10\\
11.85 2.6861397215801e-10\\
11.9 2.47081919212917e-10\\
11.95 2.27281787348629e-10\\
12 2.090833885906e-10\\
12.05 1.92310396570699e-10\\
12.1 1.76897481442698e-10\\
12.15 1.62716250951688e-10\\
12.2 1.49690944929768e-10\\
12.25 1.37678067855828e-10\\
12.3 1.2664307962449e-10\\
12.35 1.16490401374996e-10\\
12.4 1.07151986452648e-10\\
12.45 9.85626920963048e-11\\
12.5 9.06657718531197e-11\\
12.55 8.33931311873216e-11\\
12.6 7.67071918643237e-11\\
12.65 7.05581006774121e-11\\
12.7 6.49077682848621e-11\\
12.75 5.97004666923283e-11\\
12.8 5.4923382252142e-11\\
12.85 5.050692890141e-11\\
12.9 4.64572530349227e-11\\
12.95 4.27357065197336e-11\\
13 3.93063678831879e-11\\
13.05 3.61550306577429e-11\\
13.1 3.32647329657773e-11\\
13.15 3.05868121549533e-11\\
13.2 2.8131724444195e-11\\
13.25 2.58781417230698e-11\\
13.3 2.37993724728595e-11\\
13.35 2.18905238579305e-11\\
13.4 2.01343811771598e-11\\
13.45 1.85186589957916e-11\\
13.5 1.70345650961414e-11\\
13.55 1.56660968223657e-11\\
13.6 1.44082001187762e-11\\
13.65 1.32518583857675e-11\\
13.7 1.21882277560944e-11\\
13.75 1.12100778005831e-11\\
13.8 1.03099629024594e-11\\
13.85 9.48353860464722e-12\\
13.9 8.72126504424415e-12\\
13.95 8.02148546087324e-12\\
14 7.37694484742526e-12\\
14.05 6.78366614227193e-12\\
14.1 6.240577962333e-12\\
14.15 5.73745397424361e-12\\
14.2 5.27682069016715e-12\\
14.25 4.85481630082187e-12\\
14.3 4.46822167480055e-12\\
14.35 4.10445094129422e-12\\
14.4 3.77522926758309e-12\\
14.45 3.47042567866468e-12\\
14.5 3.20062938465623e-12\\
14.55 2.93497051224893e-12\\
14.6 2.69891715257636e-12\\
14.65 2.48258609465723e-12\\
14.7 2.2819117123948e-12\\
14.75 2.09823936395405e-12\\
14.8 1.9301646152946e-12\\
14.85 1.78314647424168e-12\\
14.9 1.63154654462831e-12\\
14.95 1.50077644045428e-12\\
15 1.388235474721e-12\\
15.05 1.2782568545051e-12\\
15.1 1.1665080609761e-12\\
15.15 1.07163162204797e-12\\
15.2 9.97748521633855e-13\\
15.25 9.0566791803983e-13\\
15.3 8.33257312048547e-13\\
15.35 7.69230360919686e-13\\
15.4 7.03718998771202e-13\\
15.45 6.53335322672208e-13\\
15.5 5.95727105067107e-13\\
15.55 5.47026119837514e-13\\
15.6 5.02299734598752e-13\\
15.65 4.63460227414392e-13\\
15.7 4.2454412602827e-13\\
15.75 3.91004042021296e-13\\
15.8 3.60170286311755e-13\\
15.85 3.33038614407932e-13\\
15.9 3.027375813591e-13\\
15.95 2.7859672230025e-13\\
16 2.57216478403353e-13\\
16.05 2.35021130871704e-13\\
16.1 2.19909939394191e-13\\
16.15 1.98664628466537e-13\\
16.2 1.82282876129682e-13\\
16.25 1.67273471302884e-13\\
16.3 1.53792797354141e-13\\
16.35 1.42290918949145e-13\\
16.4 1.29670090565929e-13\\
16.45 1.21518649652793e-13\\
16.5 1.09438980898751e-13\\
16.55 1.01820884313189e-13\\
16.6 9.34363782292953e-14\\
16.65 8.81360175072204e-14\\
16.7 7.87734650991539e-14\\
16.75 7.14655270944504e-14\\
16.8 7.15572143524283e-14\\
16.85 6.5977240898979e-14\\
16.9 5.53583522304689e-14\\
16.95 5.1465947676353e-14\\
17 4.74347604543919e-14\\
17.05 4.28413933491019e-14\\
17.1 3.91222968166965e-14\\
17.15 3.64929559978817e-14\\
17.2 5.28922437045651e-14\\
17.25 3.05580724021468e-14\\
17.3 3.64127069063435e-14\\
17.35 2.77870785748856e-14\\
17.4 2.35412077686529e-14\\
17.45 2.49137759179822e-14\\
17.5 2.03768004424657e-14\\
17.55 1.84578718122212e-14\\
17.6 2.00752077073861e-14\\
17.65 1.51769997188752e-14\\
17.7 1.65292748053426e-14\\
17.75 1.27883404219135e-14\\
17.8 1.75923906158672e-14\\
17.85 1.22694516513942e-14\\
17.9 1.03558978681583e-14\\
17.95 1.10250906036801e-14\\
18 8.35068645614254e-15\\
18.05 2.10993831996341e-14\\
18.1 8.02866496454523e-15\\
18.15 6.61866111782884e-15\\
18.2 5.85448749298875e-15\\
18.25 2.25397494019543e-14\\
18.3 5.23561629455928e-15\\
18.35 4.90517062327023e-15\\
18.4 4.40824118642722e-15\\
18.45 3.81922382788433e-15\\
18.5 9.90142811669029e-15\\
18.55 3.94559070411429e-15\\
18.6 3.90775018905154e-15\\
18.65 4.2504451457313e-15\\
18.7 4.53812015144238e-15\\
18.75 4.73973006339413e-15\\
18.8 2.39973127364693e-14\\
18.85 3.94964707554531e-15\\
18.9 1.76781636431558e-15\\
18.95 2.94428317976411e-15\\
19 1.5306976270935e-15\\
19.05 5.01804346067136e-15\\
19.1 5.13812892469063e-15\\
19.15 2.53289067495174e-15\\
19.2 2.5441030595867e-15\\
19.25 1.18351927861104e-15\\
19.3 2.50324766279564e-15\\
19.35 1.28918465680012e-15\\
19.4 1.5679093668252e-15\\
19.45 1.05919480861466e-15\\
19.5 1.20890705221442e-14\\
19.55 2.43324094471221e-15\\
19.6 1.07203361567385e-15\\
19.65 2.15901660530756e-15\\
19.7 4.71149543212617e-15\\
19.75 1.16165570646557e-15\\
19.8 1.04570843466459e-14\\
19.85 5.8837788190706e-16\\
19.9 6.46936542655075e-16\\
19.95 1.12889525264972e-14\\
20 5.82611656406017e-16\\
20.05 1.9210746187516e-15\\
20.1 2.6619910991578e-15\\
20.15 6.97559838736573e-15\\
20.2 1.76486008072143e-14\\
20.25 1.47969151885846e-15\\
20.3 1.28729497784779e-15\\
20.35 1.46715526845499e-15\\
20.4 2.00316989202043e-15\\
20.45 1.31973300725995e-14\\
20.5 4.66716679520145e-16\\
20.55 7.31911768712887e-16\\
20.6 1.8325168871719e-15\\
20.65 3.1772526166877e-15\\
20.7 1.84788596233514e-15\\
20.75 2.54573166715489e-15\\
20.8 5.24851925563311e-15\\
20.85 4.94392104886482e-16\\
20.9 7.09944141615932e-16\\
20.95 2.63856393125998e-14\\
21 3.66477494677313e-16\\
21.05 5.54485486029457e-18\\
};
\end{axis}
\end{tikzpicture}
	\caption{Semilogarithmic plot of the free energy versus time.}
	\label{fig.ent}
\end{figure}
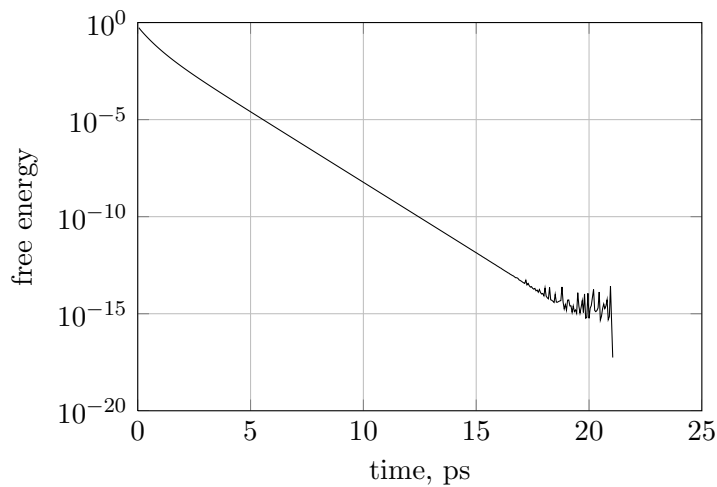

Finally, we illustrate the free energy decay in Figure \ref{fig.ent}.
In this experiment, we have set $V_D=0$ (source-drain voltage) and $V_G=0$
(source-gate voltage). It turns out that the free energy decays with an 
exponential rate. For times larger than about 18\,ps, the steady state
is almost attained, and we observe numerical oscillations caused by 
the finite machine precision.

%%%%%%%%%%%%%%%%%%%%%%%%%%%%%%%%%%%%%%%%%%%%%%%%%%%%%%%%%%%%%%%%%%%%%%%%%%%%%%%


\begin{thebibliography}{11}
\bibitem{AHPPRS14} C.~Abert, G.~Hrkac, M.~Page, D.~Praetorius, M.~Ruggeri, and
D.~Suess. Spin-polarized transport in ferromagnetic multilayers: An unconditionally
convergent FEM integrator. {\em Computers Math. Appl.} 68 (2014), 639-654.

\bibitem{ARBVHPS14} C.~Abert, M.~Ruggeri, F.~Bruckner, C.~Vogler, G.~Hrkac, 
D.~Praetorius, and D.~Suess. Self-consistent micromagnetic simulations including 
spin-diffusion effects. Preprint, 2014. {\tt arXiv:1410:6067}.

\bibitem{BaMe10} L.~Barletti and F.~M\'ehats. Quantum drift-diffusion modeling
of spin transport in nanostructures. {\em J. Math. Phys.} 51 (2010), 053304,
20 pages.

\bibitem{BCV14} M.~Bessemoulin-Chatard, C.~Chainais-Hillairet, and M.-H.~Vignal.
Study of a finite volume scheme for the drift-diffusion system. Asymptotic
behavior in the quasi-neutral regime. {\em SIAM J. Numer. Anal.} 52-4 (2014), 
1666-1691.

\bibitem{Cha11} M.~Chatard. Asymptotic behavior of the Scharfetter-Gummel scheme
for the drift-diffusion model. Proceedings of the conference ``Finite Volumes
for Complex Applications. VI. Problems and Perspectives''. {\em Springer Proc. Math.}
4 (2011), 235-243.

\bibitem{ElH14} R.~El Hajj. Diffusion models for spin transport derived from
the spinor Boltzmann equation. {\em Commun. Math. Sci.} 12 (2014), 565-592.

\bibitem{EGH00} R.~Eymard, T.~Gallou\"et, and R.~Herbin. Finite volume methods. 
In: P.~G.~Ciarlet and J.~L.~Lions (eds.). {\em Handbook of Numerical Analysis}, 
Vol.~7. North-Holland, Amsterdam (2000), 713-1020.

\bibitem{FMESZ07} J.~Fabian, A.~Matos-Abiague, C.~Ertler, P.~Stano, and 
I.~\v{Z}uti\'c. Semiconductor spintronics. {\em Acta Phys. Slovava} 57 (2007),
565-907.

\bibitem{GaGa96} H.~Gajewski and K.~G\"artner. On the discretization of van 
Roosbroeck's equations with magnetic field. {\em Z. Angew. Math. Mech.} 76
(1996), 247-264.

\bibitem{GaGl10} K.~G\"artner and A.~Glitzky. Existence of bounded steady state 
solutions to spin-polarized drift-diffusion systems. 
{\em SIAM J. Math. Anal.} 41 (2010), 2489-2513.

\bibitem{Gli08} A.~Glitzky. Analysis of a spin-polarized drift-diffusion model.
{\em Adv. Math. Sci. Appl.} 18 (2008), 401-427.

\bibitem{Gli11} A.~Glitzky. Uniform exponential decay of the free energy for
Voronoi finite volume discretized reaction-diffusion systems. {\em Math. Nachr.}
284 (2011), 2159-2174.

\bibitem{HJP04} S.~Holst, A.~J\"ungel and P.~Pietra. An adaptive mixed scheme for 
energy-transport simulations of field-effect transistors. 
{\em SIAM J. Sci. Comp.} 25 (2004), 1698-1716.

\bibitem{Ili69} A.M.~Il'in. A difference scheme for a differential equation with a 
small parameter affecting the highest derivative. {\em Mat. Zametki} 6 (1969),
237-248 (in Russian).

\bibitem{Jue09} A.~J\"ungel. {\em Transport Equations for Semiconductors}. 
Lecture Notes in Physics 773. Springer, Berlin, 2009.

\bibitem{JNS15} A.~J\"ungel, C.~Negulescu, and P.~Shpartko. Bounded weak solutions 
to a matrix drift-diffusion model for spin-coherent electron transport in 
semiconductors. To appear in {\em Math. Models Meth. Appl. Sci.}, 2015.

\bibitem{LLN08} T.~Low, M.~Lundstrom, and D.~Nikonov. Modeling of spin 
metal-oxide-semiconductor field-effect transistor: A nonequilibrium Green's
function approach with spin relaxation. {\em J. Appl. Phys.} 104 (2008), 094511, 
10 pages.

\bibitem{PSP} Y.~Pershin, S.~Saikin, and V.~Privman. Semiclassical transport models
for semiconductor spintronics. {\em Electrochem. Soc. Proc.} 2004-13 (2005), 183-205.

\bibitem{PoNe11} S.~Possanner and C.~Negulescu. Diffusion limit of a generalized
matrix Boltzmann equation for spin-polarized transport. {\em Kinetic Related Models}
4 (2011), 1159-1191.

\bibitem{ScGu69} D.I.~Scharfetter and H.K.~Gummel, Large-signal analysis of a silicon 
Read diode oscillator. {\em IEEE Trans. Electron Dev.} ED-16 (1969), 64-77.

\bibitem{Sel84} S.~Selberherr. {\em Analysis and Simulation of Semiconductor Devices.} 
Springer, Vienna, 1984.

\bibitem{ZFD02} I.~\v{Z}uti\'c, J.~Fabian, and S.~Das Sarma. Spin-polarized 
transport in inhomogeneous magnetic semiconductors: theory of 
magnetic/nonmagnetic p-n junctions. {\em Phys. Rev. Lett.} 88 (2002), 066603.

\bibitem{ZFD04} I.~\v{Z}uti\'c, J.~Fabian, and S. Das Sarma. Spintronics: 
Fundamentals and applications. {\em Rev. Modern Phys.} 76 (2004), 323-410.

\end{thebibliography}
\end{document}